\documentclass[oneside,11pt,reqno]{amsart}
\usepackage{stmaryrd}
\usepackage{bbm}
\usepackage{mathrsfs}
\usepackage{latexsym,amsxtra,xypic}
\usepackage[dvips]{graphicx}
\usepackage{dsfont}
\usepackage[all]{xy}
\usepackage{amscd,graphics}
\usepackage{amsmath,amsfonts,amsthm,amssymb}
\usepackage{slashed}
\usepackage{latexsym,amsmath}
\usepackage{graphicx,psfrag}
\usepackage{array}
\usepackage{hyperref}

\newtheorem{thm}{Theorem}[section]
\newtheorem{lem}[thm]{Lemma}
\newtheorem{cor}[thm]{Corollary}
\newtheorem{pro}[thm]{Proposition}
\newtheorem{ex}[thm]{Example}
\newtheorem{rmk}[thm]{Remark}
\newtheorem{defi}[thm]{Definition}
\newtheorem{fact}[thm]{Fact}

\title
{Pin(2)-equivariant KO-theory and intersection forms of spin four-manifolds}

\author{Jianfeng Lin}

\begin{document}

\maketitle 
\begin{abstract}
Using Seiberg-Witten Floer spectrum and Pin$(2)$-equivariant KO-theory, we prove new Furuta-type inequalities on the intersection forms of spin cobordisms between homology $3$-spheres. As an application, we give explicit constrains on the intersection forms of spin $4$-manifolds bounded by Brieskorn spheres $\pm\Sigma(2,3,6k\pm1)$. Along the way, we also give an alternative proof of Furuta's improvement of $10/8$-theorem for closed spin-$4$ manifolds.
\end{abstract}

\section{Introduction}
A natural question in 4-dimensional topology is: which nontrivial symmetric bilinear form can be realized as the intersection form of a closed, smooth, spin 4-manifold $X$. Such form should be even and unimodular. Therefore, it is indefinite by Donaldson's diagonalizability theorem \cite{Donaldson1,Donaldson2}. After changing the orientation of $X$ if necessary, we can assume that the signature $\sigma(X)$ is non-positive. Then the intersection form can be decomposed as $p(-E_{8})\oplus q\left(\begin{smallmatrix}
 0 & 1 \\
 1& 0
\end{smallmatrix}\right)$
 with  $p \geq 0, q>0$. Matsumoto's $11/8$ conjecture \cite{Matsumoto} states that $b_{2}(X)\geq \frac{11}{8}|\sigma(X)|$, which can be rephrased as $q\geq \frac{3p}{2}$. An important result is the following $10/8$ theorem of Furuta.

 \begin{thm}[Furuta \cite{Furuta}]
 Suppose $X$ is an oriented closed spin four-manifold with intersection form $p(-E_{8})\oplus q\left(\begin{smallmatrix}
 0 & 1 \\
 1& 0
\end{smallmatrix}\right)$ for $p \geq 0, q>0$. Then we have $q\geq p+1$.
 \end{thm}

Furuta's proof made use of the finite dimensional approximation of the Seiberg-Witten equations on closed 4-manifolds and Pin($2$)-equivariant $K$-theory. By doing destabilization and appealing to a result by Stolz \cite{Stolz}, Minami \cite{Minami} and Schmidt \cite{Schmidt} independently proved the following improvement:

\begin{thm}[Minami \cite{Minami}, Schmidt \cite{Schmidt}]\label{inequality for closed mfd}
 Let $X$ be a smooth, oriented, closed spin four-manifold with intersection form $p(-E_{8})\oplus q\left(\begin{smallmatrix}
 0 & 1 \\
 1& 0
\end{smallmatrix}\right)$ for $p \geq 0, q>0$. Then we have:

\begin{equation}
 q\geq \left\{
   \begin{array}{l}
 p+1,\text{ }p\equiv 0,2 \text{ mod }8  \\
 p+2,\text{ }p\equiv 4 \text{ mod }8 \\
 p+3,\text{ }p\equiv 6 \text{ mod }8. \\
   \end{array}
  \right.
  \end{equation}

\end{thm}

\begin{rmk}
$p$ is always an even integer by Rokhlin's theorem \cite{Rokhlin}.
\end{rmk}

An interesting observation is that Schmidt's calculation in \cite{Schmidt} about the Adams operations actually implies an alternative proof of the following further improvement, which was first proved by Furuta-Kametani \cite{Furuta2}. We will give the proof in Section 3.

\begin{thm}[Furuta-Kametani \cite{Furuta2}]\label{further improvement}
 Let $X$ be a smooth, oriented, closed spin four-manifold with intersection form $p(-E_{8})\oplus q\left(\begin{smallmatrix}
 0 & 1 \\
 1& 0
\end{smallmatrix}\right)$ for $p,q>0$. Then $q\geq p+3$ when $p\equiv0$ mod $8$.
\end{thm}

Another direction is to consider the intersection form of a spin $4$-manifold with given boundary. Suppose $X$ is not closed but has boundary components, which are homology three-spheres. The intersection form of $X$ is still even and unimodular but can be definite now. For the definite case, various constrains are found in \cite{Froyshov1, Froyshov2,Froyshov3,Ozvath,KMOS,Manolescu3}.

For the indefinite case, Furuta-Li \cite{Furuta-Li} and Manolescu \cite{Manolescu2} proved the following theorem independently\footnote{We give Manolescu's statement here. Furuta-Li's statement is slightly different.}.

\begin{thm}[Furuta-Li \cite{Furuta-Li}, Manolescu \cite{Manolescu2}]
To each oriented homology $3$-sphere $Y$, we can associate an invariant $\kappa(Y)\in \mathds{Z}$ with the following properties:

(i) Suppose $W$ is a smooth, spin cobordism from $Y_{0}$ to $Y_{1}$, with intersection form $p(-E_{8})\oplus q\left(\begin{smallmatrix}
 0 & 1 \\
 1& 0
\end{smallmatrix}\right)$. Then:
$$\kappa(Y_{1})+q\geq \kappa(Y_{0})+p-1.$$

(ii) Suppose $W$ is a smooth, oriented spin manifold with a single boundary $Y$, with intersection form $p(-E_{8})\oplus q\left(\begin{smallmatrix}
 0 & 1 \\
 1& 0
\end{smallmatrix}\right)$ and $q>0$. Then:
$$\kappa(Y)+q\geq p+1.$$
\end{thm}

Both Furuta-Li and Manolescu proved this theorem by considering the Pin($2$)-equivariant K-theory on the Seiberg-Witten Floer spectrum. Some new bounds can be obtained from this theorem. For example, the Brieskorn sphere $+\Sigma(2,3,12n+1)$ does not bound a spin $4$-manifold with intersection form $p(-E_{8})\oplus p\left(\begin{smallmatrix}
 0 & 1 \\
 1& 0
\end{smallmatrix}\right)$ for $p>0$.

The main purpose of this paper is to extend Theorem \ref{inequality for closed mfd} to the case of spin cobordisms and get more constrains on the intersection form of a spin $4$-manifold with boundary. Here is the first result:
\begin{thm}\label{inequality}
For any $k\in \mathds{Z}/8$, we can associate an invariant $\kappa o_{k}(Y)$ to each oriented homology sphere $Y$, with the following properties:
 \begin{itemize}
\item (1) $2\kappa o_{k}(Y)$ is an integer whose mod $2$ reduction is the Rokhlin invariant $\mu(Y)$.
\item (2) Suppose $W$ is an oriented smooth spin cobordism from $Y_{0}$ to $Y_{1}$, with intersection form $
p(-E_{8})\oplus q\left(\begin{smallmatrix}
 0 & 1 \\
 1& 0
\end{smallmatrix}\right)
$ for $p,q\geq0$. Let $p=4l+m$ for $l\in \mathds{Z}$ and $m=0,1,2,3$. Then for any $k\in \mathds{Z}/8$, we have the following inequalities:

(i) If $(\mu(Y_{0}),m)=(0,0),(0,3),(1,0),(1,1)$, then:
\begin{equation}
\kappa o_{k}(Y_{0})+2l+h(\mu(Y_{0}),m)\leq \kappa o_{k+q}(Y_{1})+\beta_{k+q}^{q}.
\end{equation}

(ii) If $(\mu(Y_{0}),m)=(0,1),(0,2),(1,2),(1,3)$, then:
\begin{equation}
\kappa o_{k+4}(Y_{0})+2l+h(\mu(Y_{0}),m)\leq \kappa o_{k+q}(Y_{1})+\beta_{k+q}^{4+q}.
\end{equation}
Here $\beta_{k}^{j}=\sum\limits_{i=0}^{j-1}\alpha_{k-i}$ where $\alpha_{i}=1$ for $i\equiv 1,2,3,5$ mod $8$ and $\alpha_{i}=0$ for $i\equiv 0,4,6,7$ mod $8$ ($\beta_{k}^{0}$ is defined to be $0$). The constants $h(\mu(Y_{0}),m)$ are listed below:

\end{itemize}
\begin{center}
    \begin{tabular}{ |c|c | c |c |c |}
    \hline
     & $m=0$ & $m=1$ & $m=2$ & $m=3$ \\ \hline
    $\mu(Y_{0})=0$ & $0$ & $5/2$ & $3$ &$3/2$ \\ \hline
    $\mu(Y_{0})=1$ & $0$ & $1/2$ & $3$ &$7/2$ \\ \hline
    \end{tabular}.
\end{center}
\end{thm}

\begin{rmk}
When $m$ is even, $\mu(Y_{0})=\mu(Y_{1})$ and $h(\mu(Y_{0}),m)$ is an integer. When $m$ is odd, $\mu(Y_{0})\neq\mu(Y_{1})$ and $h(\mu(Y_{0}),m)$ is a half-integer.
\end{rmk}

Setting $p=q=0$ in (2) of Theorem \ref{inequality}, we get:
\begin{cor}\label{cobordant invariantce}
If two homology spheres $Y_{0}, Y_{1}$ are homology cobordant to each other, then $\kappa o_{k}(Y_{0})=\kappa o_{k}(Y_{1})$ for any $k\in \mathds{Z}/8$.
\end{cor}

The definition of $\kappa o_{k}$ is similar to that of $\kappa$ \cite{Furuta-Li,Manolescu2}. Roughly, $\kappa o_{k}(Y)$ is defined as follows. Pick a metric $g$ on $Y$. By doing finite dimensional approximation to the Seiberg-Witten equations on $(Y,g)$, we get a topological space $I_{\nu}$ with an action by $G=Pin(2)$. After changing $I_{\nu}$ by suitable suspension or desuspension, we consider the following construction: The inclusion of the $S^{1}$-fixed point set $I_{\nu}^{S^{1}}$ induces a map between the equivariant KO-groups $i^{*}:\widetilde{KO}_{G}(I_{\nu})\rightarrow \widetilde{KO}_{G}(I_{\nu}^{S^{1}})$. We choose a specific reduction $\varphi:\widetilde{KO}_{G}(I_{\nu}^{S^{1}})\rightarrow \mathds{Z}$. It can be proved that the image of $\varphi\circ i^{*}$ is an ideal generated by $2^{a}\in \mathds{Z}$. We define $a$ as $\kappa o_{k}(Y)$. Different $k\in \mathds{Z}/8$ correspond to different suspensions.

In Section 8, we calculate some examples using the results in \cite{Manolescu2} about the Seiberg-Witten Floer spectrum of $\pm\Sigma(2,3,r)$.
\begin{thm}\label{kappao for Brieskorn}
(a) We have $\kappa o_{i}(S^{3})=0$ for any $i\in \mathds{Z}/8$.

(b) For a positive integer $r$ with $\text{gcd}(r,6)=1$, let $\Sigma(2,3,r)$ be the Brieskorn spheres oriented as boundaries of negative plumbings and let $-\Sigma(2,3,r)$ be the same Brieskorn spheres with the orientations reversed. Then $\kappa o_{i}(\pm\Sigma(2,3,r))$ are listed below:
\begin{center}
    \begin{tabular}{ |c|c|c |c |c |c|c|c|c|}
    \hline
     & $\kappa o_{0}$ & $\kappa o_{1}$ & $\kappa o_{2}$ & $\kappa o_{3}$ & $\kappa o_{4}$& $\kappa o_{5}$& $\kappa o_{6}$& $\kappa o_{7}$\\ \hline
    $\Sigma(2,3,12n-1)$ & $1$ & $1$ & $1$ & $0$ & $0$ & $0$ & $0$ & $0$\\ \hline
    $-\Sigma(2,3,12n-1)$ & $0$ & $0$ & $-1$ & $-1$ & $0$ & $0$ & $0$ & $0$\\ \hline
    $\Sigma(2,3,12n-5)$ & $1/2$ & $1/2$ & $1/2$ & $-1/2$ & $-1/2$ & $-1/2$ & $-1/2$ & $-1/2$\\ \hline
    $-\Sigma(2,3,12n-5)$ & $3/2$ & $3/2$ & $1/2$ & $-1/2$ & $-1/2$ & $-1/2$ & $-1/2$ & $1/2$\\ \hline
    $\Sigma(2,3,12n+1)$ & $0$ & $0$ & $0$ & $0$ & $0$ & $0$ & $0$ & $0$\\ \hline
    $-\Sigma(2,3,12n+1)$ & $0$ & $0$ & $0$ & $0$ & $0$ & $0$ & $0$ & $0$\\ \hline
    $\Sigma(2,3,12n+5)$ & $3/2$ & $3/2$ & $1/2$ & $-1/2$ & $-1/2$ & $-1/2$ & $1/2$ & $3/2$\\ \hline
    $-\Sigma(2,3,12n+5)$ & $-1/2$ & $-1/2$ & $-1/2$ & $-1/2$ & $-1/2$ & $-1/2$ & $-1/2$ & $-1/2$\\ \hline
    \end{tabular}.
\end{center}

\end{thm}

\begin{rmk}
We see that $\kappa o_{k}(-Y)\neq -\kappa o_{k}(Y)$ in general, while $\kappa o_{k}(Y\#(-Y))$ is always $0$ by Corollary \ref{cobordant invariantce}. Therefore, $\kappa o_{k}$ is not additive under connected sum.
\end{rmk}

 If we apply (2) of Theorem \ref{inequality} to the case $Y_{0}=Y_{1}=S^{3}$, the result is weaker than Theorem \ref{inequality for closed mfd}. As the case in K-theory (See \cite{Manolescu2}), we can remedy this by considering the special property of $Y_{0}\cong S^{3}$, which is called the Floer $KO_{G}$-split condition.

\begin{thm}\label{split inequality for mfd}
 Let $W$ be an oriented, smooth spin cobordism from $Y_{0}$ to $Y_{1}$, with intersection form $
p(-E_{8})\oplus q\left(\begin{smallmatrix}
 0 & 1 \\
 1& 0
\end{smallmatrix}\right)
$ and $p\geq 0,q> 0$. Suppose $Y_{0}$ is Floer $KO_{G}$-split. Let $p=4l+m$ for $l\in \mathds{Z}$ and $m=0,1,2,3$. Then we have the following inequalities:

(1)If $(\mu(Y_{0}),m)=(0,0),(0,3),(1,0),(1,1)$, then:
\begin{equation}
\kappa o_{4}(Y_{0})+2l+h(\mu(Y_{0}),m)+1\leq \kappa o_{4+q}(Y_{1})+\beta_{4+q}^{q}.
\end{equation}

(2)If $(\mu(Y_{0}),m)=(0,1),(0,2),(1,2),(1,3)$, then:
\begin{equation}
\kappa o_{4}(Y_{0})+2l+h(\mu(Y_{0}),m)+1\leq \kappa o_{q}(Y_{1})+\beta_{q}^{4+q}.
\end{equation}
Here $\beta_{*}^{*}$ and $h(\mu(Y_{0}),m)$ are the constants defined in Theorem \ref{inequality}.
\end{thm}
In particular, $S^{3}$ is Floer $KO_{G}$-split. Applying $Y_{0}=S^{3}$ to the previous theorem, we get the following useful corollary:
\begin{cor}\label{single boundary}
 Let $W$ be an oriented smooth spin $4$-manifold whose boundary is a homology sphere $Y$. Suppose the intersection form of $W$ is $
p(-E_{8})\oplus q\left(\begin{smallmatrix}
 0 & 1 \\
 1& 0
\end{smallmatrix}\right)
$ with $p\geq 0,q> 0$. Then we have the following inequalities:
\begin{itemize}
\item If $p=4l$, then $2l< \kappa o_{4+q}(Y)+\beta_{4+q}^{q}$.
\item If $p=4l+1$, then $2l+\frac{5}{2}< \kappa o_{q}(Y)+\beta_{q}^{4+q}$.
\item If $p=4l+2$, then $2l+3< \kappa o_{q}(Y)+\beta_{q}^{4+q}$.
\item If $p=4l+3$, then $2l+\frac{3}{2}< \kappa o_{4+q}(Y)+\beta_{4+q}^{q}$.
\end{itemize}
\end{cor}

\begin{rmk}
 If we set $Y=S^{3}$ in Corollary \ref{single boundary}, we will recover Theorem \ref{inequality for closed mfd}.
 However, Corollary \ref{single boundary} is not enough to prove Theorem \ref{further improvement}. In order to get the relative version of Theorem \ref{further improvement}, we have to apply similar constructions on the fixed point set of the Adams operation. This will not be done in the present paper.
\end{rmk}

 Combining the results in Theorem \ref{kappao for Brieskorn} with Corollary \ref{single boundary}, we get some new explicit bounds on the intersection forms of spin four-manifolds bounded by $\pm\Sigma(2,3,r)$. We give two of them here and refer to Section 8.2 for a complete list.

\begin{ex} We have the following conclusions:

\begin{itemize}
\item $-\Sigma(2,3,12n-1)$ does not bound a spin four-manifold with intersection form  $
p(-E_{8})\oplus (p+1)\left(\begin{smallmatrix}
 0 & 1 \\
 1& 0
\end{smallmatrix}\right)
$ for $p> 0$.
\item $-\Sigma(2,3,12n-5)$ does not bound a spin four-manifold with intersection form  $
p(-E_{8})\oplus p\left(\begin{smallmatrix}
 0 & 1 \\
 1& 0
\end{smallmatrix}\right)
$ for $p> 1$.
\end{itemize}
\end{ex}

The paper is organized as follows: In Section 2, we discuss some background material about Pin($2$)-equivariant KO-theory. In Section 3, we prove Theorem \ref{further improvement} after recalling some basic facts and properties of the Adams operations. In Section 4, we review the basic properties of the Seiberg-Witten Floer spectrum. The numerical invariant $\kappa o_{k}$ is defined in Section 5 and Theorem \ref{inequality} is proved in Section 6. In Section 7, we introduce the Floer $KO_{G}$-split condition and prove Theorem \ref{split inequality for mfd}. In Section 8, we prove Theorem \ref{kappao for Brieskorn} and use Corollary \ref{single boundary} and Theorem \ref{further improvement} to obtain new constrains on the intersection form of a spin four-manifold with given boundary.

\bigskip\noindent\textbf{Acknowledgement} Many of the constructions are parallel to those in \cite{Furuta-Li,Manolescu2} and are credited throughout. I wish to thank Ciprian Manolescu for suggesting the problem that leads to the results in this paper, and for his encouragement and enthusiasm. I am also grateful to the referee for comments on a previous version of this paper.

\section{Equivariant KO-theory}

\subsection{General Theory}
In this subsection, we review some general facts about equivariant KO-theory, mostly from \cite{Segal} and \cite{Atiyah4}. See \cite{Atiyah1}, \cite{Atiyah2} for basic facts about ordinary K-theory and KO-theory.

Let $G$ be a compact topological group and $X$ be a compact $G$-space. We denote the Grothendieck group of real  $G$-bundles over $X$ by $KO_{G}(X)$.

\begin{fact}
$KO_{G}(\text{pt})=RO(G)$. Here $RO(G)$ denotes the real representation ring of $G$. For a general $X$, $KO_{G}(X)$ is a $RO(G)$-algebra (with unit).
\end{fact}
\begin{rmk}
In this paper, we will not distinguish a representation of $G$ with its representation space.
\end{rmk}

\begin{fact}
A continuous $G$-map $f:X\rightarrow Y$ induces a map $f^{*}:KO_{G}(Y)\rightarrow KO_{G}(X)$.
\end{fact}

\begin{fact}\label{restriction}
For each subgroup $H\subseteq G$, by restricting the $G$ action to $H$, which makes a $G$-bundle into an $H$-bundle, we get a functorial restriction map $r:KO_{G}(X)\rightarrow KO_{H}(X)$.
\end{fact}
\begin{fact}
If $G$ acts freely on $X$, then the pull back map $KO(X/G)\rightarrow KO_{G}(X)$ is a ring isomorphism.
\end{fact}

\begin{fact}
For a real irreducible representation space $V$ of $G$, $\text{End}_{G}(V)$ is either $\mathds{R}$, $\mathds{C}$ or $\mathds{H}$. Let $\mathds{Z}\text{Ir}_{\mathds{R}}$, $\mathds{Z}\text{Ir}_{\mathds{C}}$ and $\mathds{Z}\text{Ir}_{\mathds{H}}$ denote the free abelian groups generated by irreducible representations of respective types and let $KSp(X)$ be the the Grothendieck group of quaternionic vector bundles over $X$. Then if $G$ acts trivially on $X$, we have:
 \begin{equation}
 KO_{G}(X)=(KO(X)\otimes\mathds{Z}\text{Ir}_{\mathds{R}})\oplus (K(X)\otimes\mathds{Z}\text{Ir}_{\mathds{C}})\oplus (KSp(X)\otimes\mathds{Z}\text{Ir}_{\mathds{H}}).
 \end{equation}
\end{fact}

Now suppose $X$ has a distinguished base point $p$ which is fixed by $G$. Then we define $\widetilde{KO}_{G}(X)$ (the reduced KO-group) to be the kernel of the map $KO_{G}(X)\rightarrow KO_{G}(p)$. For based space $X$ with trivial action, we also have:
\begin{equation}\label{trivial action}
 \widetilde{KO}_{G}(X)=(\widetilde{KO}(X)\otimes\mathds{Z}\text{Ir}_{\mathds{R}})\oplus (\widetilde{K}(X)\otimes\mathds{Z}\text{Ir}_{\mathds{C}})\oplus (\widetilde{KSp}(X)\otimes\mathds{Z}\text{Ir}_{\mathds{H}}).
 \end{equation}

The following fact is proved as Corollary 3.1.6 in \cite{Atiyah1}. (\cite{Atiyah1} only proved the complex K-theory case but the proof works without modification in the real case.)
\begin{fact}
 Suppose $X$ is a finite, based $G$-CW complex and the $G$-action is free away from the base point. Then any element in $\widetilde{KO}_{G}(X)\cong \widetilde{KO}(X/G)$ is nilpotent.
\end{fact}
Recall that the augmentation ideal $\mathfrak{a}\subset RO(G)$ is the kernel of the forgetful map $RO(G)\cong KO_{G}(pt)\rightarrow KO(pt)\cong \mathds{Z}$. Any element in $\mathfrak{a}$ defines an element in $\widetilde{KO}_{G}(X)$. By the above fact, we get:
\begin{fact}\label{nilpotent}
 Suppose $X$ is a finite, based $G$-CW complex and the $G$-action is free away from the base point. Then any element in the augmentation ideal acts on $\widetilde{KO}_{G}^{*}(X)$ nilpotently.
\end{fact}
\begin{fact}
For pointed spaces $X,Y$, there is a natural product map $\widetilde {KO}_{G}(X)\otimes\widetilde {KO}_{G}(Y)\rightarrow \widetilde{KO}_{G}(X\wedge Y)$.
\end{fact}

\begin{fact}
For pointed spaces $X,Y$, we have $\widetilde {KO}_{G}(X\vee Y)\cong \widetilde {KO}_{G}(X)\oplus \widetilde {KO}_{G}(Y)$
\end{fact}
Let $V$ be a real representation space of $G$. Denote the reduced suspension $ V^{+}\wedge X$ by $\Sigma^{V}X$. The following equivariant version of real Bott periodicity theorem was proved in \cite{Atiyah4}.

\begin{fact}
Suppose the dimension $n$ of $V$ is divisible by $8$ and $V$ is a spin representation (which means the group action $G\rightarrow SO(n)\subset End(V)$ factors through $Spin(n)$). Then we have the Bott isomorphism $\varphi_{V}:\widetilde{KO}_{G}(X)\cong\widetilde{KO}_{G}(\Sigma^{V}X)$, given by the multiplication of the Bott Class $b_{V}\in  \widetilde {KO}_{G}(V^{+})$ under the natural map $\widetilde{KO}_{G}(V^{+})\otimes \widetilde{KO}_{G}(X)\rightarrow \widetilde{KO}_{G}(\Sigma^{V}X)$. Bott isomorphism is funtorial under the pointed map $X\rightarrow X'$.
\end{fact}

\begin{fact}
Bott classes behave well under the restriction map, which means that
$i^{*}b_{V}= b_{i^{*}(V)}$. Here $i^{*}$ is the restriction map (see Fact \ref{restriction}) and $i^{*}(V)$ is the the restriction of the representation to the subgroup.
\end{fact}

\subsection{Pin(2)-equivariant KO-theory}

In this section, we will review some important facts about Pin($2$)-equivariant KO-theory. The detailed discussions can be found in \cite{Schmidt}. From now on, we assume $G\cong Pin(2)$ unless otherwise noted. Recall that the group Pin($2$) can be defined as $S^{1}\oplus jS^{1}\subset \mathds{C}\oplus j \mathds{C}=\mathds{H}$. We have:
$$RO(Pin(2))\cong \mathds{Z}[D,K,H]/(D^{2}-1,DK-K,DH-H,H^{2}-4(1+D+K)).$$

The representation space of $D$ is $\mathds{R}$ where the identity component $S^{1}\subset Pin(2)$ acts trivially and $j\in Pin(2)$ act as multiplication by $-1$.

The representation space of $K$ is $\mathds{C}\cong \mathds{R}\oplus i\mathds{R}$ where $z\in S^{1}\subset Pin(2)$ acts as multiplication by $z^{2}$ (in $\mathds{C}$) and $j$ acts as reflection along the diagonal.

The representation space of $H$ is $\mathds{H}$ where the action is given by the left multiplication of $Pin(2)\subset \mathds{H}$.

We will also write $\mathds{R}$ as the trivial one dimensional representation of $G$.

Following the notation of \cite{Schmidt}, we denote $\widetilde{KO}_{G}((kD+lH)^{+})$ by $KO_{G}(kD+lH)$ (we choose $\infty$ as the base point). Then for $k,l,m,n \in \mathds{Z}_{\geq 0}$ we have the multiplication map:
\begin{equation}\label{multiply}
KO_{G}(kD+lH)\otimes KO_{G}(mD+nH)\rightarrow KO_{G}((k+m)D+(l+n)H).
\end{equation}

In order to define this map, we need to fix the identification between $(kD\oplus lH)\oplus (mD\oplus nH) $ and $(k+m)D\oplus (l+n)H$ by sending $(x_{1}\oplus y_{1})\oplus(x_{2} \oplus y_{2})$ to $(x_{1},x_{2})\oplus (y_{1},y_{2})$. By considering the $G$-equivariant homotopy, it is not hard to see that the multiplication map is commutative when $k$ or $l$ is even. (We will prove that the multiplication map is actually commutative for any $k,l$, after we give the structure of $KO_{G}(kD+lH)$ in Theorem \ref{Pin(2)-gourp}.)

It is easy to prove (see \cite{Schmidt}) that $8D$, $H+4D$ and $2H$ are spin representations. Therefore, we can choose Bott classes $b_{8D}\in KO_{G}(8D)$, $b_{2H}\in KO_{G}(2H)$ and $b_{H+4D}\in KO_{G}(H+4D)$. Multiplication by these classes induces isomorphism $KO_{G}(kD+lH)\cong KO_{G}((k+8)D+lH)\cong KO_{G}((k+4)D+(l+1)H)\cong KO_{G}(kD+(l+2)H)$. Since the Bott classes are in the center, it doesn't matter whether we multiply on the left or on the right. Moreover, we can choose the Bott classes to be compatible with each other, which means that $b_{8D}b_{2H}=b_{H+4D}^{2}$. We will fix the choice of these Bott classes throughout this paper.

For $k,l\in \mathds{Z}$, the $RO(G)$-module $KO_{G}(kD+lH)$ is defined to be $KO_{G}((k+8a)D+(l+2b)H)$ for any $a,b\in \mathds{Z}$ which make $k+8a\geq0$ and $l+2b\geq 0$. Since the Bott Classes are chosen to be compatible, the groups defined by different choices of $a,b$ are canonically identified to each other. Again because the Bott classes are in the center, the multiplication map (\ref{multiply}) can now be extended to all $k,l,m,n \in \mathds{Z}$.

Consider the inclusion $i:7D^{+}\rightarrow 8D^{+}$. There is a unique element $\gamma(D)\in KO_{G}(-D)$ which satisfies $\gamma(D)b_{8D}=i^{*}(b_{8D})$. The map $ KO_{G}((k+1)D+lH)\stackrel{\cdot\gamma(D)}{\longrightarrow} KO_{G}(kD+lH)$ is just the map induced by the inclusion $kD\oplus lH\rightarrow (k+1)D\oplus lH$ for $k,l\geq 0$. Similarly, we can define $\gamma(H)\in KO_{G}(-H)$ and $\gamma(H+4D)=\gamma(H)\gamma(D)^{4}$. Since left multiplication and right multiplication by $\gamma(D)$ or $\gamma(H)$ just correspond to different inclusions of subspaces, which are homotopic to each other, we see that $\gamma(D)$ and $\gamma(H)$ are both in the center.

By Bott periodicity, we only have to compute $KO_{G}(lD)$ for $l=-2,-1,0,...,5$. This was done in \cite{Schmidt} and we list the result here:

\begin{thm}
[Schmidt \cite{Schmidt}]\label{Pin(2)-gourp}
As $\mathds{Z}$-modules we have the following isomorphisms:
\begin{itemize}
\item 1) $KO_{G}(pt)\cong RO(Pin(2))\cong \mathds{Z}[D,A,B]/(D^{2}-1,DA-A,DB-B,B^{2}-4 (A-2B))$, where $A=K-(1+D)$ and $B=H-2(1+D)$. \footnote{There is a typo in \cite{Schmidt}, where the relation between $A$ and $B$ is $B^{2}-2 (A-2B)$.}
\item 2) $KO_{G}(-lD)\cong\mathds{Z}\oplus\oplus_{n\geq 1}\mathds{Z}/2$ for $l=1,2$ generated by $\gamma(D)^{|l|}$ and $\gamma(D)^{|l|}A^{n}$.
\item 3) $KO_{G}(D)\cong \mathds{Z}$, generated by $\eta(D)$.
\item 4) $KO_{G}(lD)\cong \mathds{Z}\oplus\oplus_{m\geq 0}\mathds{Z}/2$ for $l=2,3.$ The generators are $\eta(D)^{2}$ and $\gamma(D)^{2}A^{m}c$ for $l=2$; $\gamma(D)\lambda(D)$ and $\gamma(D)A^{m}c$ for $l=3$.
\item 5) $KO_{G}(4D)$ is freely generated by $\lambda(D), D\lambda(D), A^{n}\lambda(D)$ and $A^{m}c$ for $m\geq0$ and $n\geq 1$.
\item 6) $KO_{G}(5D)\cong \mathds{Z}$, generated by $\eta(D)\lambda(D)$.

\end{itemize}

\end{thm}

\begin{cor}
The multiplication map (\ref{multiply}) is commutative.
\end{cor}
\begin{proof}
We just need to check $\gamma(D)$, $\eta(D)$, $\lambda(D)$, $c$ commute with each other. This is easy since $\lambda(D)$ and $c$ are in $KO_{G}(kD)$ for even $k$, while $\gamma(D)$ is in the center by our discussion before. \end{proof}

For our purpose, we don't need to know the explicit constructions of $\eta(D),\lambda(D)$ and $c$. We just need to know the following properties of them.

$\eta(D)$ is the Hurewicz image of an element $\tilde{\eta}(D)\in\pi_{G}^{0}(D)$ ($G$-equivariant stable cohomotopy group of $D^{+}$). If we forget about the $G$-action, $\tilde{\eta}(D)$ is just the Hopf map in $\pi_{1}^{\text{st}}(\text{pt})$.

For $\lambda(D)$ and $c\in KO_{G}(4D)$, by Bott periodicity and formula (\ref{trivial action}), we have isomorphisms:
$$
KO_{G}(4D)\cong KO_{G}(8D+4)\cong KO_{G}(4)$$
$$\cong (\widetilde{KO}(S^{4})\otimes\mathds{Z}\text{Ir}_{\mathds{R}})\oplus (\widetilde{K}(S^{4})\otimes\mathds{Z}\text{Ir}_{\mathds{C}})\oplus (\widetilde{KSp}(S^{4})\otimes\mathds{Z}\text{Ir}_{\mathds{H}}).
$$
(Here $4\in RO(G)$ denotes the trivial $4$-dimensional real representation. )

We can choose suitable Bott classes such that under these isomorphisms, $\lambda(D)$ corresponds to $([V_{H}]-4\mathds{R})\otimes1
\in \widetilde{KO}(S^{4})\otimes\mathds{Z}\text{Ir}_{\mathds{R}}$ and $c$ corresponds to $([V_{\mathds{H}}]-\mathds{H})\otimes H
\in \widetilde{KSp}(S^{4})\otimes\mathds{Z}\text{Ir}_{\mathds{H}}$. Here $V_{\mathds{H}}$ is the quaternion Hopf bundle over $S^{4}\cong \mathds{H}P^{2}$. $\mathds{H}$ and $\mathds{R}$ denote the trivial bundles and $1,H$ are elements in $RO(G)$.

Let $\lambda(H)$ and $c(H)$ be the image of $\lambda(D)$ and $c$ under the Bott isomorphism $KO_{G}(4D)\cong KO_{G}(8D+H)\cong KO_{G}(H)$. Then $KO_{G}(H)$ is generated by $\lambda(H)$ and $c(H)$ as $RO(G)$-algebra.

\begin{rmk}
Notice that the element $[V_{H}]\otimes H\in KSp{S^{4}}\otimes \mathds{Z}Ir_{\mathds{H}}$ is represented by the bundle $V_{H}\otimes_{\mathds{H}} H$. Hence it is a real bundle of dimension $4$ (not $16$).
\end{rmk}

For further discussions, we need to know the multiplicative structures of $KO_{G}(lD)$, which are also given in \cite{Schmidt}. We list some of them that are useful for us:
\begin{thm}
[Schmidt \cite{Schmidt}]\label{multiplicative structure} The following relations hold:
\begin{itemize}
\item 1) $H\lambda(D)=4c$, $Hc=(A+2+2D)\lambda(D)$, $Dc=c.$
\item 2) $(D+1)\gamma(D)=2A\gamma(D)=B\gamma(D)=0.$
\item 3) $(D+1)\eta(D)=A\eta(D)=B\eta(D)=0.$
\item 4) $\gamma(D)\eta(D)=1-D$, $\gamma(D)\lambda(D)=\eta(D)^{3}.$
\item 5) $\gamma(D)^{8}b_{8D}=8(1-D),\text{ }\gamma(H)^{2}b_{2H}=K-2H+D+5.$
\item 6) $\gamma(H+4D)b_{H+4D}=4(1-D).$
\item 7) $\eta(D)\lambda(D)=\gamma(D)^{3}b_{8D}$, $\eta(D)c=0.$
\item 8) $\gamma(H)\lambda(H)=4-H$ and $\gamma(H)c(H)=H-1-D-K.$

\end{itemize}
\end{thm}

\section{The Adams operations}

\subsection{Basic properties}
In this subsection, we give a quick review about the basic properties of the Adams operations. See \cite{Atiyah1} and \cite{Tom Dieck2} for more detailed discussions. Some of the calculations can be found in \cite{Schmidt} but we give them here for completeness. For simplicity and concreteness, we only deal with  $\psi^{k}:KO_{G}(X)\rightarrow KO_{G}(X) $ for an actual $G$-space $X$ and we don't do localizations (like \cite{Schmidt}).

Let $KO_{G}(X)[[t]]$ be the formal power series with coefficients in $KO_{G}(X)$. For a bundle $E$ over $X$, we define $\lambda_{t}(E)\in KO_{G}(X)[[t]]$ to be $\sum\limits_{i=0}t^{i}[\lambda^{i}(E)]$. Here $\lambda^{i}(E)$ is the $i$-th exterior power of $E$. We let $\psi^{0}(E)=\text{rank}(E)$ and define $\psi_{t}(E)=\sum\limits_{i=0}t^{i}\psi^{i}(E)\in KO_{G}(X)[[t]]$ by
\begin{equation}\label{Adams operations}
\psi_{t}(E)=\psi^{0}(E)-t\frac{d}{dt}(\text{log} \lambda_{-t}(x)).
\end{equation}
It turns out that for any $k\in \mathds{Z}_{\geq 0}$, $\psi^{k}$ extends to a well defined operation on $KO_{G}(X)$, which satisfies the following nice properties:
\begin{itemize}
\item (1) $\psi^{k}$ is functorial with respect to continuous maps $f:X\rightarrow X'$.
\item (2) $\psi^{k}$ maps $\widetilde{KO}_{G}(X)$ to $\widetilde{KO}_{G}(X)$.
\item (3) For all $x,y\in KO_{G}(X)$, $\psi^{k}(x+y)=\psi^{k}(x)+\psi^{k}(y)$ and $\psi^{k}(xy)=\psi^{k}(x)\psi^{k}(y)$.
\item (4) If $x$ is a line bundle, then $\psi^{k}(x)=x^{k}$.
\end{itemize}
The effect of the Adams operations on the Bott classes can be described by the Bott cannibalistic class. Given a spin $G$-bundle $E$ over $X$ with rank $n\equiv0$ mod $8$, the Bott cannibalistic class $\theta^{\text{or}}_{k}(E)\in RO(G)$ is defined by the equation:
\begin{equation}\label{defi of cannibalistic}
\psi^{k}(b_{E})=\theta^{\text{or}}_{k}(E)\cdot b_{E} \text{ for }k>1.
\end{equation} When $k$ is odd, this can be explicitly written as (see \cite{Tom Dieck2}):\footnote{There is a typo in 3.10.4 \cite{Tom Dieck2}.}
\begin{equation}\label{Bott Cannibalistic Class}
\theta^{\text{or}}_{k}(E)=k^{n/2}\prod_{u\in J}\lambda_{-u}(E)(1-u)^{-n}.
\end{equation}
Here $J$ is a set of $k$-th unit roots $u\neq1$ such that $J$ contains exactly one element from each pair $\{u,u^{-1}\}$.
Notice that we can define $\theta^{\text{or}}_{k}(E)$ for any real bundle $E$ of even dimension using formula (\ref{Bott Cannibalistic Class}). It can be shown that: $$\theta^{\text{or}}_{k}(E+F)=\theta^{\text{or}}_{k}(E)\theta^{\text{or}}_{k}(F).$$

Now let's specialize to the case $k=3$. By formula ({\ref{Adams operations}}), it is easy to check that $\psi^{3}(x)=x^{3}-3\lambda^{2}(x)x+3\lambda^{3}(x)$. We want to calculate the action of $\psi^{3}$ on $RO(G)$. Since the $G$-action on $H$ preserves the orientation, we have $\lambda^{3}(H)=\lambda^{1}(H)=H$. Using complexification, it is easy to show $\lambda^{2}(H)=K+D+3$. Also, we have $\lambda^{2}(K)=D$. Therefore, we get\footnote{There is a typo in \cite{Schmidt}, where $\psi^{3}(H)=HK-K$.}:
$$\psi^{3}(D)=D,\text{ }\psi^{3}(H)=HK-H,\text{ }\psi^{3}(K)=K^{3}-3K,$$
$$\psi^{3}(A)=A^{3}+6A^{2}+9A,\text{ }\text{ }\psi^{3}(B)=AB+B+4A.$$

Also, applying formula (\ref{Bott Cannibalistic Class}), we get:
$$\theta^{\text{or}}_{3}(2)=3,\text{ }\theta^{\text{or}}_{3}(2D)=1+2D,\text{ }\theta^{\text{or}}_{3}(H)=A+B+4D+5.$$
\subsection{Proof of Theorem \ref{further improvement}}
The central part of the proof is the following proposition:
\begin{pro}\label{stable map between spheres}
For any integers $r,a,b\geq 0$ and $l>0$, there does not exist $G$-equivariant map
$$f:(r\mathds{R}+aD+(4l+b)H)^{+}\rightarrow (r\mathds{R}+(a+8l+2)D+bH)^{+}$$
which induces homotopy equivalence on the $G$-fixed point set.
\end{pro}
\begin{proof}
Suppose there exists such a map $f$. After suspension by copies of $\mathds{R},D$ and $H$, we can assume $a=8l'+6$, $r=8d$ and $b=2k$. Let $V_{1}=8d\mathds{R}+2kH+8(l+l'+1)D$ and $V_{2}=8d\mathds{R}+(4l+2k)H+(8l'+8)D$. Let $b_{V_{1}}$ and $b_{V_{2}}$ be the Bott classes of $V_{1}$ and $V_{2}$, respectively. Consider the element $x=f^{*}(b_{V_{1}})$. By the Bott isomorphism and (2) of Theorem \ref{Pin(2)-gourp}, we can write $x$ as $b_{V_{2}}\gamma(D)^{2}\alpha$ for some $\alpha\in RO(G)$. Moreover, we can assume $\alpha=p+Ah(A)$ for some integer $p$ and some polynomial $h(A)$ whose coefficients are either $0$ or $1$.

\bigskip\textbf{Claim}: $p$ is even and $h=0$.

This is essentially a special case of Proposition 5.21 in \cite{Schmidt} for $\mathcal{KO}(4l,8l+2)$.\footnote{ There is an error in \cite{Schmidt} for $\mathcal{KO}(c,d)$ when $4c-d\equiv-3$ mod $8$, but we will not consider this case here.}

By formula (\ref{defi of cannibalistic}), we have: $\psi^{3}(b_{V_{1}})=\theta_{3}^{\text{or}}(V_{1})\cdot b_{V_{1}}$, which implies:
\begin{equation}\label{method 1}\psi^{3}(x)=f^{*}(\psi^{3}(b_{V_{1}}))=\theta_{3}^{\text{or}}(V_{1})\cdot x.\end{equation}
Notice that $x=i^{*}(b_{V_{2}}\cdot\alpha)$ where $i: (8d\mathds{R}+(4l+2k)H+(8l'+6)D)^{+}\rightarrow V_{2}^{+}$ is the standard inclusion. By formula (\ref{defi of cannibalistic}), we have:
\begin{equation}\label{method 2}\psi^{3}(x)=i^{*}(\psi^{3}(b_{V_{2}}\cdot\alpha))=\theta_{3}^{\text{or}}(V_{2})b_{V_{2}}\psi^{3}(\alpha)\cdot \gamma(D)^{2}.\end{equation}

Comparing equation (\ref{method 1}) and equation (\ref{method 2}), we get:
\begin{equation}\label{method 1=2}
(\theta_{3}^{\text{or}}(V_{2})\psi^{3}(\alpha)-\theta_{3}^{\text{or}}(V_{1})\alpha)\gamma(D)^{2}=0
\end{equation}

We can calculate:
$$
\theta_{3}^{\text{or}}(V_{1})=3^{4d}(1+2D)^{4l+4l'+4}(A+B+4D+5)^{2k},
$$
$$
\theta_{3}^{\text{or}}(V_{2})=3^{4d}(1+2D)^{4l'+4}(A+B+4D+5)^{2k+4l}.
$$

Notice that $2A\gamma(D)=B\gamma(D)=(1+D)\gamma(D)=0$, we can simplify equation (\ref{method 1=2}) as:
\begin{equation}\label{simplification 1}
3^{4d}((A+1)^{2k}\alpha-(A+1)^{4l+2k}\psi^{3}(\alpha))\cdot\gamma(D)^{2}=0.
\end{equation}
Since $\alpha=p+Ah(A)$, we have $\psi^{3}(\alpha)=p+(A^{3}+6A^{2}+9A)h(A^{3}+6A^{2}+9A)$.
  Using the relation $2A\gamma(D)=0$, we can further simplify equation (\ref{simplification 1}) and get:
\begin{equation}
3^{4d}\cdot g(A)\cdot \gamma(D)^{2}=0
\end{equation}
Here $g(A)=(A+1)^{2k}(p+Ah(A))-(A+1)^{2k+4l}(p+(A^{3}+A)h(A^{3}+A)).$

By (2) of Theorem \ref{Pin(2)-gourp}, we see that if we expand $g(A)$ as a polynomial in $A$, the degree-$0$ coefficient should be $0$ and all other coefficients should be even. By our assumption, the coefficients of $h$ are either $0$ or $1$. Checking the leading coefficient of $g(A)$, it is easy to see that $h=0$ and $g(A)=p((A+1)^{2k}-(A+1)^{2k+4l})$. This implies that $p$ is even. The claim is proved.

Now consider the commutative diagram:
\begin{equation}
\xymatrix{
\widetilde{KO}_{G}(V_{1}^{+}) \ar[d]^{\cdot\gamma(H)^{2k}\gamma(D)^{8l+8l'+8}} \ar[r]^-{f^{*}} & \widetilde{KO}_{G}((8d\mathds{R}+(8l'+6)D+(4l+2k)H)^{+})\ar[d]^{\cdot\gamma(H)^{4l+2k}\gamma(D)^{8l'+6}}\\
\widetilde{KO}_{G}((8d\mathds{R})^{+}) \ar[r]^-{\cong} & \widetilde{KO}_{G}((8d\mathds{R})^{+}).}
\end{equation}
The vertical maps are given by the inclusions of subspaces. The bottom map is an isomorphism because $f$ induces a homotopy equivalence on the $G$-fixed point set. Any automorphism on $\widetilde{KO}_{G}((8d\mathds{R})^{+})$ is given by the multiplication of a unit $\tilde{u}\in RO(G)$. Therefore, we obtain :
\begin{equation}\tilde{u}\cdot b_{V_{1}}\cdot \gamma(H)^{2k}\gamma(D)^{8l+8l'+8}= x\cdot \gamma(H)^{4l+2k}\gamma(D)^{8l'+6}=b_{V_{2}}\cdot\gamma(D)^{8l'+8}\gamma(H)^{4l+2k}\cdot p\end{equation}

Applying the relations in Theorem \ref{multiplicative structure}, we simplify this as :
\begin{equation}
(K-2H+D+5)^{2l+k}(8(1-D))^{l'+1}\cdot p=(K-2H+D+5)^{k}(8(1-D))^{l+l'+1}\cdot \tilde{u}.
\end{equation}

Now consider the ring homomorphism $\varphi_{0}:RO(G)\rightarrow \mathds{Z}$ defined by $\varphi_{0}(D)=-1,\varphi_{0}(A)=\varphi_{0}(B)=0$. Notice that $\varphi_{0}(\tilde{u})=\pm1$ since $\tilde{u}$ is a unit. We get $p=\pm1$, which is a contradiction. This finishes the proof of Proposition \ref{stable map between spheres}.\end{proof}

 Now suppose $X$ is a closed, oriented, smooth spin four-manifold with intersection form $p(-E_{8})\oplus q\left(\begin{smallmatrix}
 0 & 1 \\
 1& 0
\end{smallmatrix}\right)$ for $p=8l>0$ and $q<p+3$. After doing surgery on loops and connect sum copies of $S^{2}\times S^{2}$, we can assume $b_{1}(W)=0$ and $q=8l+2$. As shown in \cite{Furuta}, by doing finite dimensional approximation of the Seiberg-Witten equations on $W$, we get an $G$-equivariant map:
$$f:(aD+(4l+b)H)^{+}\rightarrow ((a+8l+2)D+bH)^{+}\text { for some } a,b>0.$$
 Moreover, $f$ induces a homotopy equivalence on the $G$-fixed point set. This is a contradiction to Proposition \ref{stable map between spheres}. Therefore, Theorem \ref{further improvement} is proved.

\section{Pin(2)-equivariant Seiberg-Witten Floer theory}
In \cite{Manolescu1}, \cite{Manolescu2} and \cite{Manolescu3}, Manolescu constructed a Pin(2)-equivariant spectrum class $S(Y,\mathfrak{s})$ for each rational homology sphere $Y$ with a spin structure $\mathfrak{s}$. We will not repeat the constructions here but just collect some useful properties. See \cite{Manolescu1}, \cite{Manolescu2} and \cite{Manolescu3} for the explicit constructions.

\begin{defi}
Let $s\in \mathds{Z}_{\geq 0}$. A space of type SWF (at level s) is a pointed, finite G-CW complex $X$ with the following properties:
\begin{itemize}
\item(a) The $S^{1}$-fixed point set $X^{S^{1}}$ is $G$-homotopy equivalent to the sphere $(sD)^{+}$. We define $\text{lev}(X)$ to be $s$.
\item(b) The action of $G$ is free on the complement $X-X^{S^{1}}$.
\end{itemize}
\end{defi}

\begin{defi}
Let $X,X'$ be two spaces of type SWF at level $k$ and $k'$ respectively. A pointed $G$-map $f:X\rightarrow X'$ is called admissible if $f$ preserves the base point and satisfies one of the following two conditions:
\begin{itemize}
\item (1) $k<k'$ and the induced map on the $G$-fixed point set $f^{G}:X^{G}\rightarrow X'^{G}$ is a homotopy equivalence.
\item (2) $k=k'$ and the induced map on the $S^{1}$-fixed point set $f^{S^{1}}:X^{S^{1}}\rightarrow X'^{S^{1}}$ is a homotopy equivalence.
\end{itemize}

\end{defi}

Now consider the set of triples $(X,a,b)$ where $X$ is a space of type SWF and $a\in \mathds{Z}, b\in \mathds{Q}$.

\begin{defi}\label{stable equivalence}
We say that $(X,a,b)$ is stable equivalent to $(X',a',b')$ if $b-b'\in \mathds{Z}$ and for some $M,N,r>0$, there exists a $G$-homotopy equivalence:
$$\Sigma^{r\mathds{R}}\Sigma^{(M-a)D}\Sigma^{(N-b)H}X\cong \Sigma^{r\mathds{R}}\Sigma^{(M-a')D}\Sigma^{(N-b')H}X'.$$
(Here $\mathds{R}$ denotes the trivial representation of $G$.)
\end{defi}
\begin{rmk}
In \cite{Manolescu2}, Manolescu worked with stable even equivalence, which requires $X$ to be a space of type SWF at even level.
\end{rmk}

This triple can be thought of the ``formal de-suspension'' of $X$ with $a$ copies of $D$ and $b$ copies of $H$. We denote $\mathfrak{C}$ to be the set of stable equivalence classes of triples $(X,a,b)$. Informally, we call an element in $\mathfrak{C}$ a spectrum class.

\begin{defi}
For a spectrum class $S= [(X,a,b)]\in \mathfrak{C}$, we let $$\text{lev}(S)=\text{lev}(X)-a.$$
\end{defi}

\begin{rmk}
By considering the $S^{1}$-fixed point set, we see that two spaces of type SWF at different levels are not $G$-homotopic to each other. Using this fact, it is easy to prove that $lev(S)$ is a well defined quality.
\end{rmk}

For $r\in \mathds{Z}$ and $s\in \mathds{Q}$, we can define the formal suspension $\Sigma^{rD+sH}:\mathfrak{C}\rightarrow \mathfrak{C}$ by sending $[(X,a,b)]$ to $[(X,a-r,b-s)]$. It's easy to check that this is a well defined operation on the set $\mathfrak{C}$.

Now suppose $Y$ is an oriented rational homology three-sphere with a metric $g$ and a spin structure $\mathfrak{s}$. Let $\mathbb{S}$ be the associated spinor bundle. We consider the global Coulomb splice:
 $$
 V=i\text{ker} d^{*}\oplus \Gamma(\mathbb{S})\subset i\Omega^{1}(Y)\oplus \Gamma(\mathbb{S}).
 $$
 Using the quaternionic structure on $\mathbb{S}$, we can define a natural action of $G$ on $V$: $e^{i\theta}\in G$ takes $(\alpha,\phi)$ to $(e^{i\theta}\alpha,\phi)$ and $j\in G$ takes $(\alpha,\phi)$ to $(-\alpha,j\phi)$.

 Now we consider the self-adjoint first order elliptic operator $l:V\rightarrow V$ defined by $l(\alpha,\phi)=(*d\alpha,\slashed{D}\phi)$ where $\slashed{D}$ is the Dirac operator \footnote{Since $Y$ is a rational homology sphere, there is a unique flat spin-connection on $\mathbb{S}$, we choose it as the base connection and use it to define $\slashed{D}$.}. For any $\tau<\nu$, let $V^{\tau}_{\nu}$ be the subspace spanned by the eigenvectors of $l$ with eigenvalues in the interval $(\tau,\nu]$. Then $V^{\tau}_{\nu}$ is a finite dimensional $G$-representation space which is isomorphic to $kD\oplus lH$. We denote $k$ by $\text{dim}_{\mathds{R}}V(D)^{\nu}_{\tau}$ and $l$ by $\text{dim}_{\mathds{H}}V(H)^{\nu}_{\tau}$.

 We pick $-\nu<<0<<\nu$. By considering the equivariant Conley index of the gradient flow of $CSD|_{V^{\nu}_{-\nu}}$ (see \cite{Manolescu1} and \cite{Manolescu2}), we get a $G$-space $I_{\nu}$ of type SWF at level $\text{dim}_{\mathds{R}}V(D)^{0}_{-\nu}$.

Next, we need to recall the definition of $n(Y,\mathfrak{s},g)$. Choose a compact smooth spin four-manifold $N$ with $\partial N=Y$. Let $\text{ind}_{\mathds{C}}\slashed{D}(N)$ be the index of Dirac operator on $N$ (with Atiyah-Patodi-Singer boundary conditions). We can define:
\begin{equation}
    n(Y,\mathfrak{s},g):=\text{ind}_{\mathds{C}}\slashed{D}(N)+\frac{\sigma(N)}{8}.
\end{equation}

\begin{rmk}
It can be proved that this definition does not depend on the choice of $N$. For a rational homology sphere $Y$, we have $n(Y,\mathfrak{s},g)\in \frac{1}{8}\mathds{Z}$. When $Y$ is an integral homology sphere, $n(Y,\mathfrak{s},g)$ is an integer and has the same parity as the Rokhlin invariant $\mu(Y)$.
\end{rmk}

We can consider the following element in $\mathfrak{C}$:\footnote{Our convention is different from \cite{Manolescu1} and \cite{Manolescu2} , where the second component in the triple denotes the complex dimension of the $G$-representation.}
\begin{equation}
S(Y,\mathfrak{s}):=[(I_{\nu},\text{dim}_{\mathds{R}}V(D)^{0}_{-\nu},\text{dim}_{\mathds{H}}V(H)^{0}_{-\nu}+\frac{1}{2}n(Y,\mathfrak{s},g))].
\end{equation}

Notice that the level of $S(Y,\mathfrak{s})$ is always $0$.

\begin{thm}
[Manolescu \cite{Manolescu1},\cite{Manolescu2}]\label{topological invariance} The element $S(Y,\mathfrak{s})\in \mathfrak{C}$ is independent of the metric $g$, the cut-off $\nu$ and the other choices in the construction. Thus  $S(Y,\mathfrak{s})$ is an invariant of the pair $(Y,\mathfrak{s})$.
\end{thm}

\begin{rmk}
In this paper, since we only use the numerical invariants, we don't need to make $\mathfrak{C}$ a category and $S(Y,\mathfrak{s})$ a functor. Therefore, we don't define $S(Y,\mathfrak{s})$ as natural spectrum invariant. See Section 3.4 of \cite{Manolescu3} for a discussion about naturality.
\end{rmk}

Suppose $W$ is a smooth spin cobordism between rational homology three spheres $Y_{0}$ and $Y_{1}$ with $b_{1}(W)=0$. Further, we assume $W$ is equipped with a matric $g$ and a spin structure $\mathfrak{t}$ such that $g|_{Y_{i}}=g_{i}$ and $\mathfrak{t}|_{Y_{i}}=\mathfrak{s}_{i}$.

The following theorem is important for our constructions:

\begin{thm}
[Manolescu \cite{Manolescu1},\cite{Manolescu2}]\label{cobordism map} By doing finite dimensional approximation for the Seiberg-Witten equations on $W$, we obtain an admissible map:
\begin{equation}
f:\Sigma^{a_{0}D}\Sigma^{b_{0}H}(I_{0})_{\nu}\rightarrow \Sigma^{a_{1}D}\Sigma^{b_{1}H}(I_{1})_{\nu}.
\end{equation}

Here, $(I_{0})_{\nu}$ and $(I_{1})_{\nu}$ are the Conley indices for the approximated Seiberg-Witten flow. Let $V_{i}$ denotes the Coulomb slice on $Y_{i}$, for $i=0,1$. The differences in the suspension indices are:
\begin{equation}
a_{0}-a_{1}=\text{dim}_{\mathds{R}}V_{1}(D)^{0}_{-\nu}-\text{dim}_{\mathds{R}}V_{0}(D)^{0}_{-\nu}-b_{2}^{+}(W)
\end{equation}
and
\begin{equation}
b_{0}-b_{1}=\text{dim}_{\mathds{H}}V_{1}(H)^{0}_{-\nu}-\text{dim}_{\mathds{H}}V_{0}(H)^{0}_{-\nu}+\frac{1}{2}n(Y_{1},\mathfrak{s}_{1},g_{1})-\frac{1}{2}n(Y_{0},\mathfrak{s}_{0},g_{0})-\frac{\sigma(W)}{16}.
\end{equation}
\end{thm}

\section{Numerical Invariants}

Let $Y$ be a rational homology sphere and $\mathfrak{s}$ be a spin structure on $Y$. In the previous section, we defined an invariant $S(Y,s)\in\mathfrak{C}$. In this section, we will extract a set of numerical invariants $\kappa o_{i}(Y,s)$ from $S(Y,s)$, for $i\in \mathds{Z}/8$.

\begin{defi}\label{reduction}
For $l=-2,-1,0,...,5$, we define the group homomorphisms $\varphi_{l}:KO(lD)\rightarrow \mathds{Z}$ as following (see Theorem \ref{Pin(2)-gourp}):

\begin{itemize}
\item 1) $\text {For } l=0$, $\varphi_{l}(D)=-1$ and $\varphi_{l}(A)=\varphi_{l}(B)=0$, then extend $\varphi_{l}$ by the multiplicative structure on $RO(G)$.
\item 2) $\text {For } l=-1,-2$, $\varphi_{l}(\gamma(D)^{|l|})=1$ and $\varphi_{l}(\gamma(D)^{|l|}A^{n})=0$ for $n\geq1$.
\item 3) $\text {For } l=1$, $\varphi_{l}(\eta(D))=1$.
\item 4) $\text {For } l=2$, $\varphi_{l}(\eta(D)^{2})=1$ and $\varphi_{l}(\gamma(D)^{2}A^{m}c)=0$.
\item 5) $\text {For } l=3$, $\varphi_{l}(\gamma(D)\lambda(D))=1$ and $\varphi_{l}(\gamma(D)A^{m}c)=0$.
\item 6) $\text {For } l=4$, $\varphi_{l}(\lambda(D))=1$, $\varphi_{l}(D\lambda(D))=-1$, and $\varphi_{l}(A^{n}\lambda(D))=\varphi_{l}(A^{m}c)=0$.
\item 7) $\text {For } l=5$, $\varphi_{l}(\eta(D)\lambda(D))=1$.

\end{itemize}
For the other $l\in \mathds{Z}$, we use the Bott isomorphism to identify $KO(lD)$ with $KO((l-8k)D)$ for $-2\leq l-8k\leq 5$ and apply the above definition.
\end{defi}

\begin{lem}\label{multiplicative of phi}
For any $a\in KO_{G}(pt)$ and $b\in KO_{G}(kD)$, we have $\varphi_{0}(a)\varphi_{k}(b)=\varphi_{k}(a\cdot b)$.
\end{lem}
\begin{proof}
This is a straightforward calculation using Theorem \ref{Pin(2)-gourp} and Theorem \ref{multiplicative structure}.
\end{proof}

\begin{rmk}
$\varphi_{0}$ is just taking the trace of $j\in Pin(2)$. While the other $\varphi_{l}$ are defined such that the torsion elements are killed and Lemma \ref{multiplicative of phi} holds.
\end{rmk}

We consider the map $\tau: D^{+}\rightarrow D^{+}$ which maps $x$ to $-x$. By suspension with copies of $D$, we get an admissible involution $\tau: (kD)^{+}\rightarrow (kD)^{+}$ for $k>0$.

The following lemma is a straightforward corollary of the equivariant Hopf theorem (see \cite{Tom Dieck}).

\begin{lem}\label{homeomorphism}
 When $0\leq k<l$, any admissible map $f:(kD)^{+}\rightarrow (lD)^{+}$ is $G$-homotopic to the standard inclusion.
For $0\leq k=l$, any admissible map $f:(kD)^{+}\rightarrow (kD)^{+}$ is either homotopic to $\tau$ or to the identity map, depending on $\text{deg}(f)$.
\end{lem}

 $\tau$ induces the involution $\tau^{*}: KO_{G}(kD)\rightarrow KO_{G}(kD)$. For $k,l>0$ and any $a\in KO_{G}(kD), b\in KO_{G}(lD)$, the following equality is easy to check by Lemma \ref{homeomorphism}:
\begin{equation}\label{tau}
\tau^{*}(a)\cdot b=a\cdot \tau^{*}(b)=\tau^{*}(a\cdot b) \text{ and } \tau^{*}(a)\cdot \tau^{*}(b)= a\cdot b.
\end{equation}

Using this fact, we can define $\tau^{*}:KO_{G}(kD)\rightarrow KO_{G}(kD)$ for any $k\in \mathds{Z}$ by identifying $KO_{G}(kD)$ with $KO_{G}(k'D)$ for any $0<k'\equiv k \text{ mod }8$ using Bott periodicity. Moreover, formula (\ref{tau}) now holds for all $k,l\in \mathds{Z}$.

Now consider the element $u\in RO(G)$ defined by $\tau^{*}(b_{8D})=u\cdot b_{8D}$. Then for $l\in\mathds{Z}$ and any element $\alpha\in KO_{G}(lD)$, we have $\tau^{*}(\alpha)\cdot b_{8D}=\alpha\cdot \tau^{*}(b_{8D})=(u\alpha)\cdot b_{8D}$, which implies $\tau^{*}(\alpha)=u\alpha$.

\begin{lem}\label{reflective invariance} We have the following properties about $\tau^{*}$ and $u$:
\begin{itemize}
\item (1) $\tau^{*}$ acts as identity on $KO_{G}(lD)$ for $l\neq 0,4$ mod $8$.
\item (2) $u$ is a unit with $\varphi_{0}(u)=1$.
\item (3) $\varphi_{l}\circ \tau^{*}=\varphi_{l}$ for any $l\in\mathds{Z}$.
\end{itemize}

\end{lem}
\begin{proof}

(1) We have $\gamma(D)b_{8D}=i^{*}(b_{8D})$ where $i^{*}$ is the inclusion $(7D)^{+}\rightarrow (8D)^{+}$. Therefore, we get $\tau^{*}(\gamma(D)b_{8D})=(\tau\circ i)^{*}(b_{8D})$. By Lemma \ref{homeomorphism}, $\tau\circ i$ is $G$-homotopic to $i$, thus $\tau^{*}(\gamma(D)b_{8D})=i^{*}(b_{8D})=\gamma(D)b_{8D}$, which implies that $\tau^{*}(\gamma(D))=\gamma(D)$.

Since $\tau^{*}$ induces an involution on $KO_{G}(D)\cong \mathds{Z}$, we have $\tau^{*}(\eta(D))=\pm \eta(D)$. But since $\tau^{*}(\eta(D))\cdot \gamma(D)=\eta(D)\cdot \tau^{*}(\gamma(D))=\eta(D)\gamma(D)=1-D\neq -\eta(D)\gamma(D)$, we get $\tau^{*}(\eta(D))=\eta(D)$.

By formula (\ref{tau}), $\tau^{*}(a)=a$ implies $\tau^{*}(ab)=ab$ for any $a,b$. Therefore we see that $\tau^{*}$ acts as the identity map on $KO_{G}(kD)$ for $k\neq 0,4$ mod $8$.

(2) $u^{2}=1$ because $\tau^{2}=\text{id}$. Since $u \cdot (1-D)=\tau^{*} (1-D)=\tau^{*} (\gamma(D)\cdot \eta(D))=\gamma(D)\cdot \eta(D)=1-D$, we see that $(u-1)(1-D)=0$. We get $\varphi_{0}(u)=1$ by Lemma \ref{multiplicative of phi}.

(3) is straightforward from (2) and Lemma \ref{multiplicative of phi}.
\end{proof}

Now suppose $X$ is a space of type SWF at level $l$. A choice of $G$-homotopy equivalence $X^{S^{1}}\cong (lD)^{+}$ gives us an inclusion map $i:(lD)^{+}\rightarrow X$, which we call a trivialization. A trivialization induces the map $i^{*}: \widetilde {KO}_{G}(X)\rightarrow KO_{G}(lD)$. Consider the map $\varphi_{l}\circ i^{*}: \widetilde {KO}_{G}(X)\rightarrow \mathds{Z}$.

\begin{pro}\label{inva}
 The submodule $\text{Im}(i^{*})$ and the map $\varphi_{l}\circ i^{*}$ are both independent of the choice of the trivialization. Moreover, we have $\text{Im}(\varphi_{l}\circ i^{*})=(2^{k})$ for some $k\in \mathds{Z}_{\geq 0}$.
\end{pro}

\begin{proof}
 By Lemma \ref{homeomorphism}, there are two possible trivializations $i$  and $i\circ \tau$. We have $\text{Im}(i\circ \tau)^{*}=\tau^{*}(\text{Im}i^{*})=u\cdot \text{Im}(i^{*})$. Since $u$ is a unit, the multiplication by $u$ does not change the submodule $\text{Im}(i^{*})$. Moreover, we have $\varphi_{l}\circ (i\circ \tau)^{*}= \varphi_{l}\circ \tau^{*}\circ i^{*}=\varphi_{l}\circ i^{*} $  by (3) of Lemma \ref{reflective invariance}.

For the second statement, we consider the exact sequence:
$$...\rightarrow \widetilde{KO}_{G}(X)\stackrel{i^{*}}{\longrightarrow} KO_{G}(lD) \stackrel{\delta}{\longrightarrow} \widetilde{KO}_{G}^{1}(X/X^{S^{1}})\rightarrow ...$$
Since the $G$ action is free away from the basepoint and $(1-D)\in RO(G)$ is in the augmentation ideal, $(1-D)$ acts on $\widetilde{KO}_{G}^{1}(X/X^{S^{1}})$ nilpotently by Fact \ref{nilpotent}. Therefore, we can find $m\gg 0$ such that $(1-D)^{m}KO_{G}(lD)\subset \text{ker}(\delta)=\text{Im}(i^{*})$. It follows that $2^{m}\in \text{Im}(\varphi_{l}\circ i^{*})$ and $\text{Im}(\varphi_{l}\circ i^{*})=(2^{k})$ for some $0\leq k\leq m$.\end{proof}

Proposition \ref{inva} justify the following definition:

\begin{defi}
For a $G$-space $X$ of type SWF at level $l$, we define $\mathcal{J}(X)$ to be the image of $i^{*}$ for any trivialization $i$ and let $\kappa o(X)$ be the integer $k$ such that $\varphi_{l}(\mathcal{J}(X))=(2^{k})$.
\end{defi}

Let's study the property of $\mathcal{J}(X)$ and $\kappa o(X)$. First recall that we defined the constants $\beta_{k}^{0}=0$ and $\beta_{k}^{j}=\sum\limits_{i=0}^{j-1}\alpha_{k-i}$ for $j\geq 1$, where $\alpha_{i}=1$ for $i\equiv 1,2,3,5$ mod $8$ and $\alpha_{i}=0$ for $i\equiv 0,4,6,7$ mod $8$. It's easy to see that $\beta_{j}^{k}= \beta_{j'}^{k}$ for $j\equiv j'(\text{mod}8)$. The integers $\beta_{j}^{k}$ are important because of the following proposition:

\begin{pro}\label{induced phi}
For integers $0\leq j\leq k$ and an admissible map $i:((k-j)D)^{+}\rightarrow (kD)^{+}$, we have the following commutative diagram, where the map $m_{k}^{j}:\mathds{Z}\rightarrow \mathds{Z}$ is the multiplication of $2^{\beta_{k}^{j}}$.
\begin{equation}
\xymatrix{
KO_{G}(kD) \ar[d]^{\varphi_{k}} \ar[r]^-{i^{*}} & KO_{G}((k-j)D)\ar[d]^{\varphi_{k-j}}\\
\mathds{Z} \ar[r]^{m_{k}^{j}} &\mathds{Z}}
\end{equation}
\end{pro}

\begin{proof}
The case $j=0$ follows from Lemma \ref{reflective invariance}. When $j>0$, by Lemma \ref{homeomorphism}, the map $i$ is $G$-homotopic to the standard inclusion. Because of the associativity of $i^{*}$ and $m_{l}^{k}$, we only need to prove the case $j=1$. In this case, the map $i^{*}$ is just the multiplication by $\gamma(D)$ and $m_{k}^{1}$ is the multiplication by $2^{\alpha_{k}}$. Since both $\varphi_{k}$ and $i^{*}$ are compatible with Bott isomorphism, we only need to check the case $k=1,2,...,8$. This can be proved by straightforward calculations using Definition \ref{reduction}, Theorem \ref{multiplicative structure} and Theorem \ref{Pin(2)-gourp}.
 \end{proof}

The following proposition studies the behavior of $\mathcal{J}(X)$ and $\kappa o(X)$ under the Bott isomorphism:

\begin{pro} \label{periodicity of kappao}
Let $X$ be a space of type SWF at level $k$. We have the following:
\begin{itemize}
\item (1) $\mathcal{J}(X)\cdot b_{8D}=\mathcal{J}(\Sigma^{8D}X)$ and $\kappa o(\Sigma^{8D}X)=\kappa o( X)$.
\item (2) $\mathcal{J}(X)\cdot (K-2H+D+5)=\mathcal{J}(\Sigma^{2H}X)$ and $\kappa o(\Sigma^{2H}X)=\kappa o(X)+2$.
\item (3) $\kappa o(\Sigma^{H+4D}X)=\kappa o(X)+3-\beta_{k+4}^{4}$.
\end{itemize}
\end{pro}

\begin{proof}
(1) Since $(\Sigma^{8D}X)^{S^{1}}=\Sigma^{8D}(X^{S^{1}})$, statement (1) follows from the functoriality of the Bott isomorphism.

(2) We have the commutative diagram induced by the inclusions of subspaces:
\begin{equation}\label{suspeension by 2H}
\xymatrix{
\widetilde{KO}_{G}(\Sigma^{2H}X) \ar[d] \ar[r] & \widetilde {KO}_{G}(X)\ar[d]\\
\widetilde{KO}_{G}((\Sigma^{2H}X)^{S^{1}}) \ar[r]^-{\cong} &\widetilde{KO}_{G}(X^{S^{1}}).}
\end{equation}

Since $(\Sigma^{2H}X)^{S^{1}}=\Sigma^{2H}(X^{S^{1}})$, the map in the bottom row is the identity. If we identify $\widetilde {KO}_{G}(\Sigma^{2H}X)$ with $\widetilde {KO}_{G}(X)$ using the Bott isomorphism, then the top horizontal map is the multiplication by $\gamma(H)^{2}b_{2H}=K-2H+D+5$ (by Theorem \ref{multiplicative structure}). This implies $\mathcal{J}(\Sigma^{2H}X)=(K-2H+D+5)\mathcal{J}(X)$. We also have $\kappa o(\Sigma^{2H}X)=\kappa o(X)+2$ since $\varphi_{0}(K-2H+D+5)=4$.

(3) Again, by inclusions of subspaces, we have:
$$
\xymatrix{
\widetilde{KO}_{G}(\Sigma^{H+4D}X) \ar[d] \ar[r] & \widetilde {KO}_{G}(X)\ar[d]\\
KO_{G}((\Sigma^{H+4D}X)^{S^{1}}) \ar[r]^-{\cdot \gamma(D)^{4}} &KO_{G}(X^{S^{1}}).}
$$

Since $(\Sigma^{H+4D}X)^{S^{1}}\cong \Sigma^{4D}(X^{S^{1}})$, the bottom horizontal map is the multiplication by $\gamma(D)^{4}$. If we identify $\widetilde{KO}_{G}(\Sigma^{H+4D}X)$ with $\widetilde {KO}_{G}(X)$ using the Bott isomorphism, the top horizontal map is the multiplication by $\gamma(H+4D)b_{H+4D}=4(1-D)$ (by Theorem \ref{multiplicative structure}). Therefore, under appropriate trivializations, we see that the maps $i_{1}^{*}: \widetilde {KO}_{G}(X)\cong \widetilde{KO}_{G}(\Sigma^{H+4D}X)\rightarrow KO_{G}((k+4)D)$ and $i_{2}^{*}:\widetilde {KO}_{G}(X)\rightarrow KO_{G}(kD)$ are related by $\gamma(D)^{4}\cdot i_{1}^{*}(x)=4(1-D)\cdot i_{2}^{*}(x)$. Since $\varphi_{0}(4(1-D))=8$, statement (3) follows from Proposition \ref{induced phi} (for $j=4$) and Lemma \ref{multiplicative of phi}.\end{proof}

We have the following proposition, which is the analogue of Lemma 3.8 in \cite{Manolescu2}.
\begin{pro}\label{suspension by trivial}
Let $X_{1}$ and $X_{2}$ be spaces of type SWF. Suppose there is a based $G$-equivariant homotopy equivalence $f$ from  $\Sigma^{r\mathds{R}}X_{1}$ to $\Sigma^{r\mathds{R}}X_{2}$, for some $r\geq0$. Then we have $\mathcal{J}(X_{1})=\mathcal{J}(X_{2})$ and $\kappa o(X_{1})=\kappa o (X_{2})$.\end{pro}
\begin{proof}
The proof in \cite{Manolescu2} works with some modifications. Suppose $X_{1}$, $X_{2}$ are both at level $k$. By (1) of Proposition \ref{periodicity of kappao}, we can replace $X_{i}$ by $\Sigma^{8D}X_{i}$ and assume $k>1$. Also, we can suspend some more copies of $\mathds{R}$ and assume that $8|r$. Choose trivilizations $i_{1},i_{2}$ of $X_{1}$ and $X_{2}$, respectively. They give homotopy equivalences $(r\mathds{R}+kD)^{+}\cong (\Sigma^{r\mathds{R}}X_{1})^{S^{1}}$ and $(r\mathds{R}+kD)^{+}\cong(\Sigma^{r\mathds{R}}X_{2})^{S^{1}}$. Composing them with $f^{S^{1}}:(\Sigma^{r\mathds{R}}X_{1})^{S^{1}}\rightarrow (\Sigma^{r\mathds{R}}X_{2})^{S^{1}}$, we get the  equivariant homotopy equivalence $h: (r\mathds{R}+kD)^{+}\rightarrow (r\mathds{R}+kD)^{+}$. Since $k>1$, by equivariant Hopf theorem, $h$ is based homotopic to $\tau_{1}\wedge \tau_{2}$. The map $\tau_{1}:(r\mathds{R})^{+}\rightarrow (r\mathds{R})^{+}$ is either identity or a map with degree $-1$. Therefore, $\tau_{1}^{*}(b_{r\mathds{R}})=a\cdot b_{r\mathds{R}}$ where $b_{r\mathds{R}}$ is the Bott class and $a\in RO(G)$ is a unit. Also, $\tau_{2}:(kD)^{+}\rightarrow (kD)^{+}$ is either identity or the map $\tau$ we defined before. Therefore, $\tau_{2}^{*}(x)$ is either $x$ or $ux$ (see Lemma \ref{reflective invariance}). We have shown that the map $h^{*}:\widetilde{KO}_{G}((r\mathds{R}+kD)^{+})\rightarrow \widetilde{KO}_{G}((r\mathds{R}+kD)^{+})$ is just multiplication by some unit in $RO(G)$, which does not change any submodule.

Now consider the following commutative diagram:
$$
\xymatrix{
\widetilde{KO}_{G}(X_{2}) \ar[d]^{i_{2}^{*}} \ar[r]^-{\cong} & \widetilde{ KO}_{G}(\Sigma^{r\mathds{R}}X_{2})\ar[d]^{(\Sigma^{r\mathds{R}}i_{2})^{*}}\ar[r]^-{f^{*}}&\ar[d]^{(\Sigma^{r\mathds{R}}i_{1})^{*}}\widetilde{ KO}_{G}(\Sigma^{r\mathds{R}}X_{1})\ar[r]^{\cong}&\widetilde{KO}_{G}(X_{1})\ar[d]^{i_{1}^{*}} \\
KO_{G}(kD) \ar[r]^-{\cong} &\widetilde{KO}_{G}((r\mathds{R}+kD)^{+})\ar[r]^-{h^{*}}&\widetilde{KO}_{G}((r\mathds{R}+kD)^{+})\ar[r]^-{\cong}&KO_{G}(kD).}
$$
In each row, the first map is a Bott isomorphism and the third map is the inverse to a Bott isomorphism. We see that $b_{r\mathds{R}}\cdot\text{Im}(i_{2}^{*})=h^{*}(b_{r\mathds{R}}\cdot\text{Im}(i_{2}^{*}))=b_{r\mathds{R}}\cdot\text{Im}(i_{1}^{*})$.
Therefore, we have $\text{Im}(i_{1}^{*})=\text{Im}(i_{2}^{*})$, which implies $\kappa o (X_{1})=\kappa o(X_{2})$.\end{proof}

\begin{defi}\label{kappao for spectrum}
For a spectrum class $S=[(X,a,b)]\in \mathfrak{C}$, we let
 \begin{equation}
 \kappa o(S)=\kappa o(\Sigma^{(8M-a)D}\Sigma^{(2N-b')H}X)-2N-s
 \end{equation}
 for any $M,N,b'\in \mathds{Z}$ and $s\in [0,1)$ making $8M-a\geq0, 2N-b'\geq0$ and $b=b'+s$.
\end{defi}

\begin{pro}
$\kappa o(S)$ is well defined.
\end{pro}

\begin{proof}
By (1) and (2) of Proposition \ref{periodicity of kappao}, it's easy to prove that the righthand side of formula (\ref{kappao for spectrum}) is independent of the choice of $M,N$. By choosing $M,N\gg 0$, we see that changing the representative of $S$ from $(X,a,b)$ to $(\Sigma^{D}X,a+1,b)$ or $(\Sigma^{H}X,a,b+1)$ does not change the value of $\kappa o(S)$. By Definition \ref{stable equivalence} and Proposition \ref{suspension by trivial}, we proved that $\kappa o(S)$ does not change when we change the representative of the spectrum class.
\end{proof}
By definition of the suspension of a spectrum class and Proposition \ref{periodicity of kappao}, it is easy to prove:
\begin{pro}\label{periodicity of kappao for spectrum}
For any spectrum class $S\in \mathfrak{C}$ at level $k$, we have:
\begin{itemize}
\item $\kappa o(\Sigma^{8D}S)=\kappa o(S)$.
\item $\kappa o(\Sigma^{2H}S)=\kappa o(S)+2$.
\item $\kappa o(\Sigma^{H+4D}S)=\kappa o(S)+3-\beta_{k+4}^{4}$.
\end{itemize}

\end{pro}

With these discussions, we can now define the invariants for three manifolds.
\begin{defi}
For an oriented rational homology sphere $Y$ and a spin structure $\mathfrak{s}$ on $Y$, we define
$\kappa o_{i}(Y,\mathfrak{s})=\kappa o(\Sigma^{iD}S(Y,\mathfrak{s}))$ for any $i\in \mathds{Z}_{\geq0}$. Then $\kappa o_{i}(Y,\mathfrak{s})=\kappa o_{i+8}(Y,\mathfrak{s})$, which allow us to define $\kappa o_{i}(Y,\mathfrak{s})$ for $i\in \mathds{Z}/8$.
\end{defi}

\section{Proof of Theorem \ref{inequality}}
In this section, we will prove Theorem \ref{inequality}.

Let $X_{0},X_{1}$ be be two spaces of type SWF at level $k_{0}$ and $k_{1}$, respectively. Suppose there is an admissible map $f:X_{0}\rightarrow X_{1}$ (which implies $k_{0}\leq k_{1}$). By Lemma \ref{induced phi}, we can choose suitable trivializations such that the following diagram commutes.
 $$
\xymatrix{
\widetilde{KO}_{G}(X_{1}) \ar[d]^{i_{1}^{*}} \ar[r]^{f^{*}} & \widetilde{KO}_{G}(X_{0})\ar[d]^{i_{0}^{*}}\\
KO_{G}(k_{1}D) \ar[r]^-{(f^{S^{1}})^{*}}\ar[d]^{\varphi_{k_{1}}} &KO_{G}(k_{0}D)\ar[d]^{\varphi_{k_{0}}}\\
\mathds{Z}\ar[r]^{m_{k_{1}}^{k_{1}-k_{0}}}& \mathds{Z}}
$$
Therefore, we get $m^{k_{1}-k_{0}}_{k_{1}}(\text{Im}(\varphi_{k_{1}}\circ i_{1}^{*}))\subset \text{Im}(\varphi_{k_{0}}\circ i_{0}^{*})$. This implies that $(2^{\kappa o(X_{1})+\beta_{k_{1}}^{k_{1}-k_{0}}})\subset (2^{\kappa o(X_{0})})\subset \mathds{Z}$. Therefore, we get the following proposition:

\begin{pro}\label{inequality for space}
Let $X_{0},X_{1}$ be two spaces of type SWF at level $k_{0}$ and $k_{1}$, respectively. Suppose there is an admissible map $f:X_{0}\rightarrow X_{1}$. Then we have:
\begin{equation}
\kappa o(X_{0})\leq\kappa o(X_{1})+\beta_{k_{1}}^{k_{1}-k_{0}}
.\end{equation}
\end{pro}

Next we generalize the above inequality to the spectrum classes:
\begin{defi}
Let $S_{0},S_{1}\in \mathfrak{C}$ be two spectrum classes. We call $S_{0}$ dominates $S_{1}$ if we can find representatives $S_{i}=[(X_{i},a,b)]$ for $i=1,2$ and an admissible map $f$ from $X_{0}$ to $X_{1}$.
\end{defi}

\begin{pro}\label{inequality for spectrum}
Let $S_{0},S_{1}\in \mathfrak{C}$ be two spectrum classes at level $k_{0}$ and $k_{1}$ respectively. Suppose $S_{0}$ dominates $S_{1}$, then we have:
\begin{equation}
\kappa o(S_{0})\leq\kappa o(S_{1})+\beta_{k_{1}}^{k_{1}-k_{0}}.
\end{equation}
\end{pro}
\begin{proof}
Since an admissible map $f:X_{0}\rightarrow X_{1}$ gives an admissible map $\Sigma^{aH+bD}f:\Sigma^{aH+bD}X_{0}\rightarrow \Sigma^{aH+bD}X_{1}$ for any $a,b\in \mathds{Z}_{\geq 0}$. This proposition is a straightforward corollary of Proposition \ref{inequality for space} and Definition \ref{kappao for spectrum}.
\end{proof}
By considering the natural inclusion $X\rightarrow \Sigma^{D}X$, it is easy to see that $S$ always dominates $\Sigma^{D}S$. Therefore ,we get the following corollary, which will be useful in Section 8.

\begin{cor}\label{relation}
For any spectrum class $S\in \mathfrak{C}$ at level $k$. We have:
$$\kappa o(S)\leq\kappa o(\Sigma^{D}S)+\alpha_{k+1}.$$
\end{cor}

Now let $Y_{0}, Y_{1}$ be two rational homology three-spheres and $\mathfrak{s}_{i}$ be spin structures on them respectively. Suppose $(W,\mathfrak{s})$ is a smooth oriented spin cobordism from $(Y_{0},\mathfrak{s}_{0})$ to $(Y_{1},\mathfrak{s}_{1})$. After doing surgery along loops in $W$, we can assume $b_{1}(W)=0$ without loss of generality. Then by Theorem \ref{cobordism map}, we see that $\Sigma^{-\frac{\sigma(W)}{16}H}S(Y_{0},\mathfrak{s}_{0})$ dominates $\Sigma^{b_{2}^{+}(W)D}S(Y_{1},\mathfrak{s}_{1})$. We can do suspensions and prove $\Sigma^{-\frac{\sigma(W)}{16}H}(\Sigma^{kD}S(Y_{0},\mathfrak{s}_{0}))$ dominates $\Sigma^{(b_{2}^{+}(W)+k)D}S(Y_{1},\mathfrak{s}_{1})$ for any $k\in \mathds{Z}$. Applying Proposition \ref{inequality for spectrum}, we get:

\begin{thm}\label{inequality for QHS}
Suppose $(W,\mathfrak{s})$ is a smooth, oriented spin cobordism from $(Y_{0},\mathfrak{s}_{0})$ to $(Y_{1},\mathfrak{s}_{1})$. Then for any $k\in \mathds{Z}$, we have the inequality:
\begin{equation}
\kappa o_{k+b_{2}^{+}(W)}(Y_{1},\mathfrak{s}_{1})+\beta_{k+b_{2}^{+}(W)}^{b_{2}^{+}(W)}\geq \kappa o(\Sigma^{-\frac{\sigma(W)}{16}H}(\Sigma^{kD}S(Y_{0},\mathfrak{s}_{0}))).
\end{equation}
\end{thm}

In general, $\kappa o(\Sigma^{-\frac{\sigma(W)}{16}H}(\Sigma^{kD}S(Y_{0},\mathfrak{s}_{0})))$ can be expressed by $\kappa o_{k}(Y_{0},\mathfrak{s}_{0})$ or $\kappa o_{k+4}(Y_{0},\mathfrak{s}_{0})$, but the explicit formula is messy. For simplicity, we now focus on the integral homology sphere case.

\begin{rmk}
Suppose $Y$ is an oriented integral homology three-sphere. There is a unique spin structure $\mathfrak{s}$ on $Y$ and we simply write $S(Y,\mathfrak{s})$ and $\kappa o_{i}(Y,\mathfrak{s})$ as $S(Y)$ and $\kappa o_{i}(Y)$, respectively.
\end{rmk}

Suppose both $Y_{i}$ are integral homology spheres, then the intersection form of $W$ is a unimodular, even form. Let's assume that the intersection from can be decomposed as:
$$
p(-E_{8})\oplus q\left(\begin{smallmatrix}
 0 & 1 \\
 1& 0
\end{smallmatrix}\right)  \text{ for }p,q\geq0.
$$

In this case, we have $\frac{\sigma(W)}{16}=-\frac{p}{2}$ and $b_{2}^{+}(W)=q$. Recall that the spectrum class invariant $S(Y_{0})$ is defined by $[(I_{\nu},\text{dim}_{\mathds{R}}V(D)^{0}_{-\nu},\text{dim}_{\mathds{H}}V(H)^{0}_{-\nu}+\frac{1}{2}n(Y_{0},\mathfrak{s},g))]$.
 The third component of this triple may be an integer or a half integer, depending on the Rokhlin invariant $\mu(Y_{0})$.
\begin{pro}
Let $Y_{0}$ be an integral homology three sphere and $p\in\mathds{Z}_{\geq 0}$. Then we have the following relations.

(1) Suppose $\mu(Y_{0})=0\in \mathds{Z}_{2}$.
\begin{itemize}
\item For $p=4l$, we have $\kappa o(\Sigma^{\frac{p}{2}H}(\Sigma^{kD}S(Y_{0})))=\kappa o_{k}(Y_{0})+2l$.
\item For $p=4l+1$, we have $\kappa o(\Sigma^{\frac{p}{2}H}(\Sigma^{kD}S(Y_{0})))=\kappa o_{k+4}(Y_{0})+\frac{5}{2}+2l-\beta_{k}^{4}.$
\item For $p=4l+2$, we have $\kappa o(\Sigma^{\frac{p}{2}H}(\Sigma^{kD}S(Y_{0})))=\kappa o_{k+4}(Y_{0})+3+2l-\beta_{k}^{4}.$
\item For $p=4l+3$, we have $\kappa o(\Sigma^{\frac{p}{2}H}(\Sigma^{kD}S(Y_{0})))=\kappa o_{k}(Y_{0})+2l+\frac{3}{2}.$
\end{itemize}

(2) Suppose $\mu(Y_{0})=1\in \mathds{Z}_{2}$.
\begin{itemize}
\item For $p=4l$, we have $\kappa o(\Sigma^{\frac{p}{2}H}(\Sigma^{kD}S(Y_{0})))=\kappa o_{k}(Y_{0})+2l.$
\item For $p=4l+1$, we have $\kappa o(\Sigma^{\frac{p}{2}H}(\Sigma^{kD}S(Y_{0})))=\kappa o_{k}(Y_{0})+2l+\frac{1}{2}.$
\item For $p=4l+2$, we have $\kappa o(\Sigma^{\frac{p}{2}H}(\Sigma^{kD}S(Y_{0})))=\kappa o_{k+4}(Y_{0})+3+2l-\beta_{k}^{4}.$
\item For $p=4l+3$, we have $\kappa o(\Sigma^{\frac{p}{2}H}(\Sigma^{kD}S(Y_{0})))=\kappa o_{k+4}(Y_{0})+\frac{7}{2}+2l-\beta_{k}^{4}.$
\end{itemize}
\end{pro}
\begin{proof} Let's denote $(I_{\nu},\text{dim}_{\mathds{R}}V(D)^{0}_{-\nu},\text{dim}_{\mathds{H}}V(H)^{0}_{-\nu}+\frac{1}{2}n(Y_{0},\mathfrak{s},g))$ by $(X,a,b)$.

 For $\mu(Y_{0})=0$ and $p=4l$, we have $ b\in \mathds{Z}$. Take $M,N\gg0$ and let $N'=N+l$. Then by Definition \ref{kappao for spectrum}, we have:
 \begin{equation}
 \begin{split}
 \kappa o(\Sigma^{\frac{p}{2}H}(\Sigma^{kD}S(Y_{0}))=\kappa o(\Sigma^{(8M+k-a)D}\Sigma^{(2N+2l-b)H}X)-2N\\
 =\kappa o(\Sigma^{(8M+k-a)D}\Sigma^{(2N'-b)H}X)-2N'+2l=\kappa o_{k}(Y)+2l.
  \end{split}\end{equation}

For $p=4l+1$, take $M,N\gg 0$ and let $N'=N+l$. Then we have:
\begin{equation}
 \begin{split}\kappa o(\Sigma^{\frac{p}{2}H}(\Sigma^{kD}S(Y_{0}))=\kappa o(\Sigma^{(8M+k-a)D}\Sigma^{(2N+2l+1-b)H}X)-2N-\frac{1}{2}\\
=\kappa o (\Sigma^{H}(\Sigma^{kD}(X,a,b)))+2l-\frac{1}{2}=\kappa o_{k+4}(Y_{0})+\frac{5}{2}+2l-\beta_{k}^{4}.
\end{split}
\end{equation}
The other cases can be proved similarly.
\end{proof}

Now combining the above proposition and Theorem \ref{inequality for QHS}, we proved Theorem \ref{inequality}.

\section{$KO_{G}$-Split condition}

Now consider the space $X=(8kD+(2l+1)H)^{+}$ for $k,l\in \mathds{Z}_{\geq 0}$. We have the map induced by the inclusion:
$$i^{*}:\widetilde {KO}_{G}(X)\rightarrow KO_{G}(8kD).$$
By Theorem \ref{Pin(2)-gourp}, we see that $KO_{G}(8kD+(2l+1)H)$ is generated by $(b_{2H})^{l}(b_{8D})^{k}\lambda(H)$ and $(b_{2H})^{l}(b_{8D})^{k}c(H)$ as $RO(G)$-module and the map $i^{*}$ is multiplication by $\gamma(H)^{2l+1}$. Using Proposition \ref{multiplicative structure}, we get:
 \begin{equation}\begin{split}
 i^{*}((b_{2H})^{l}(b_{8D})^{k}\lambda(H))=(2+A-2D-2B)^{l}(2-2D-B)\cdot(b_{8D})^{k},\\
i^{*}((b_{2H})^{l}(b_{8D})^{k}c(H))=(A-2B)^{l}(B-A)\cdot(b_{8D})^{k}. \end{split}\end{equation}

The above discussion motivates the following definition:
\begin{defi}\label{even KO-split}
   Let $X$ be a space of type SWF at level $8k$. $X$ is called even  $KO_{G}$-split if $\mathcal{J}(X)$ is the submodule generated by $(2+A-2D-2B)^{l}(2-2D-B)\cdot (b_{8D})^{k}$ and $(A-2B)^{l}(B-A)\cdot (b_{8D})^{k}$ for some $l\in \mathds{Z}_{\geq 0}$.
\end{defi}

Next, we consider the space $X=((8k+4)D+2lH)^{+}$. The map:
$$
i^{*}:\widetilde {KO}_{G}(X)\rightarrow KO_{G}((8k+4)D)
$$
is just multiplication of $\gamma(H)^{2l}$. We know $\widetilde{KO}_{G}(X)=KO_{G}((8k+4)D)\cdot (b_{2H})^{l}$ by the Bott isomorphism. Since $\gamma(H)^{2l}(b_{2H})^{l}=(K-2H+D+5)^{l}=(A+2D+6-2H)^{l}$ (see Theorem \ref{multiplicative structure}), we have $\text{Im}(i^{*})=(A+2D+6-2H)^{l}\cdot KO_{G}((8k+4)D)\subset KO_{G}((8k+4)D)$. This motivates the following definition:
\begin{defi}\label{odd KO-split}
   Let $X$ be a space of type SWF at level $8k+4$. $X$ is called odd $KO_{G}$-split if $\mathcal{J}(X)=(A+2D+6-2H)^{l}\cdot KO_{G}((8k+4)D)$  for some $l\in \mathds{Z}_{\geq 0}$.
\end{defi}

$KO_{G}$-split spaces are special because of the following proposition (compare Proposition \ref{inequality for space}).

\begin{pro}\label{split inequality for space}
Let $X_{0},X_{1}$ be two spaces of type SWF at level $k_{0},k_{1}$ respectively and $f$ be an admissible map from $X_{0}$ to $X_{1}$. Suppose $k_{0}<k_{1}$ and $X_{0}$ is odd or even $KO_{G}$-split (which implies that $k_{0}\equiv 0 \text{ or } 4 \text{ mod } 8$). Then we have:
\begin{equation}
\kappa o(X_{0})<\kappa o(X_{1})+\beta_{k_{1}}^{k_{1}-k_{0}}.
\end{equation}
\end{pro}

Before proving this proposition, we need to make a digression into the general properties of $KO_{G}(4D)$ and $RO(G)$.

\begin{lem}\label{facts about KO4} The following properties holds:
\begin{itemize}
\item(1) Any element in $RO(G)$ can be uniquely written as $bD+f(A)+Bg(A)$ for some polynomials $f,g$ and integer $b$.
\item (2) Any element in $RO(G)$ can be uniquely written as $bD+f(A)+Hg(A)$ for some polynomials $f,g$ and integer $b$.
\item (3) Any element in $KO_{G}(4D)$ can be uniquely written as $bD\lambda(D)+f(A)\lambda(D)+g(A)c$ for some polynomials $f,g$ and integer $b$.
\item (4) The map $RO(G)\rightarrow KO_{G}(4D)$ defined by multiplication of $\lambda(D)$ is injective.
\item (5) An element $\omega= bD\lambda(D)+f(A)\lambda(D)+g(A)c$ belongs to $RO(G)\lambda(D)$ if and only if $4|g(A)$. Moreover, if $(A+2D+6-2H)^{l}\omega\in RO(G)\cdot\lambda(D)$ for some $l$, then $\omega\in RO(G)\cdot\lambda(D)$.
\item (6) Suppose $(A-2B)^{l}h(A,B)=0\in RO(G)$ for some two-variable polynomial $h$ in $A,B$. Then we have $h(A,B)=0$ in $RO(G)$.
\item (7) Suppose $f(D)=h(A,B)$ for some $2$-variable polynomial $h$ without degree-$0$ term and some polynomial $f$. Then $h(A,B)=0$.
\end{itemize}
\end{lem}
\begin{proof}
(1),(2),(3),(4) can be proved by straightforward calculation using Theorem \ref{Pin(2)-gourp}. The first statement of (5) is the corollary of (2),(3) and the relation $H\lambda(D)=4c$. Let's prove the second statement of (5). We have $Hc=(1+D+K)\lambda(D)$ and $(2D+6)c=8c=2H\lambda(D)$. Therefore, $(A+2D+6-2H)^{l}\omega\in RO(G)\lambda(D)$ implies $A^{l}\omega\in RO(G)\lambda(D)$. It follows that $4|A^{l}g(A)$, which implies $4|g(A)$ and $\omega\in RO(G)\lambda(D)$.

For (6), we can assume that $h(A,B)=f(A)+Bg(A)$ for some polynomials $f,g$. Consider the map $\psi:RO(G)\rightarrow \mathds{Q}[x]$ defined by $\psi(D)=1,\psi(B)=x$ and $\psi(A)=\frac{x^{2}}{4}+2x$. Then $0=\psi((A-2B)^{l}(f(A)+Bg(A)))=(\frac {x^{2}}{4})^{l}(f(\frac{x^{2}}{4}+2x)+xg(\frac{x^{2}}{4}+2x))$, which implies $0=f(\frac{x^{2}}{4}+2x)+xg(\frac{x^{2}}{4}+2x)$. Considering the leading term in $x$, we see that $f(x)=g(x)=0$.

For (7), we can simplify $h(A,B)$ as $Ag_{1}(A)+Bg_{2}(A)$ for some polynomials $g_{1},g_{2}$ by the relation $B^{2}-4(A-2B)=0$. Then the conclusion follows from (1).
\end{proof}
\begin{lem}\label{lowest degree, odd}
Suppose $a(1-D)\lambda(D)\in (A+2D+6-2H)^{l}KO_{G}(4D)$ for some $a\in \mathds{Z}$ and $l\in \mathds{Z}_{\geq0}$. Then we have $2^{2l+1}|\varphi_{4}(a(1-D)\lambda(D))$.
\end{lem}

\begin{proof}
 Since $\varphi_{4}(a(1-D)\lambda(D))=2a$, the conclusion is trivial when $l=0$. Now suppose $l>0$. Let $a(1-D)\lambda(D)=(A+2D+6-2H)^{l}\cdot \omega$ for some $\omega\in KO_{G}(4D)$. By (5) of Lemma \ref{facts about KO4}, we see that $\omega\in RO(G)\lambda(D)$. Write $\omega$ as $(bD+f(A)+Bg(A))\lambda(D)$. By (4) of Lemma \ref{facts about KO4}, we get $a(1-D)=(A-2B-2D+2)^{l}(bD+f(A)+Bg(A))$. Using the relation $(1-D)A=(1-D)B=0$, we can simplify this equality as $a(1-D)-(f(0)+bD)(2-2D)^{l}=(A-2B)^{l}(b+f(A)+Bg(A))$. By (7) of Lemma \ref{facts about KO4}, we get that  $(A-2B)^{l}(b+f(A)+Bg(A))=0\in RO(G)$. By $(6)$ of Lemma \ref{facts about KO4}, we have $b+f(A)+Bg(A)=0$. This implies that $\omega=b(D-1)\lambda(D)$ and $\varphi_{4}(a(1-D)\lambda(D))=-2^{2l+1}b$ for some $b\in \mathds{Z}$.
 \end{proof}

\begin{lem}\label{lowest degree, even}
Suppose $a(1-D)$ is in the ideal of $RO(G)$ generated by $(2+A-2D-2B)^{l}(2-2D-B)$ and $(A-2B)^{l}(B-A)$ for some  $l\in \mathds{Z}_{\geq 0}$. Then we have $2^{2l+3}|\varphi_{0}(a(1-D))$.
\end{lem}
\begin{proof}
 We assume $l>0$ first. By (1) of Lemma \ref{facts about KO4} and the relation $A(1-D)=B(1-D)=0$, we have can express $a(1-D)$ as:\begin{equation}\begin{split}(2-2D-B)&(2-2D+A-2B)^{l}(b(1-D)+f_{1}(A)+Bg_{1}(A))\\&+(A-2B)^{l}(B-A)(f_{2}(A)+Bg_{2}(A))\end{split}\end{equation} for some integer $b$ and polynomials $f_{1},f_{2},g_{1},g_{2}$.

As in the proof of Lemma \ref{lowest degree, odd}, we can simplify this formula and use (7) of Lemma \ref{facts about KO4} to get:\begin{equation}
-B(A-2B)^{l}(f_{1}(A)+Bg_{1}(A))+(A-2B)^{l}(B-A)(f_{2}(A)+Bg_{2}(A))=0\in RO(G).
\end{equation}
We have $-B(f_{1}(A)+Bg_{1}(A))+(B-A)(f_{2}(A)+Bg_{2}(A))=0$ by $(6)$ of Lemma \ref{facts about KO4}. Simplifying this, we obtain:
\begin{equation}\label{identity for A,B}
-4Ag_{1}(A)-Af_{2}(A)+4Ag_{2}(A)+B(-f_{1}(A)+f_{2}(A)+8g_{1}(A)-Ag_{2}(A)-8g_{2}(A))=0.\end{equation}
This implies $-4Ag_{1}(A)-Af_{2}(A)+4Ag_{2}(A)=0$ and $-f_{1}(A)+8g_{1}(A)+f_{2}(A)-Ag_{2}(A)-8g_{2}(A)=0$. Considering the degree-$1$ term of the first identity, we get $4|f_{2}(0)$. Also, we have $8|-f_{1}(0)+f_{2}(0)$ by checking the degree-$0$ term of the second identity. Therefore, we have $4|f_{1}(0)$, which implies $\varphi_{0}(a(1-D))=2^{2l+2}(2b+f_{1}(0))$ can be divided by $2^{2l+3}$.

The case $l=0$ is similar. We also get the identity (\ref{identity for A,B}).
\end{proof}
\begin{proof}[Proof of Proposition \ref{split inequality for space}:]
Consider the commutative diagram:
$$
\xymatrix{
\widetilde{KO}_{G}(X_{1}) \ar[d]^{i_{1}^{*}} \ar[r]^{f^{*}} & \widetilde{KO}_{G}(X_{0})\ar[d]^{i_{0}^{*}}\\
KO_{G}(k_{1}D) \ar[r]^-{(f^{S^{1}})^{*}}\ar[d]^{\varphi_{k_{1}}} &KO_{G}(k_{0}D)\ar[d]^{\varphi_{k_{0}}}\\
\mathds{Z}\ar[r]^{m_{k_{1}}^{k_{1}-k_{0}}}& \mathds{Z}.}
$$

(1) Suppose $X_{0}$ is odd $KO_{G}$-split. Then $k_{0}=8k+4$ for some integer $k$ and $KO_{G}(k_{0}D)=KO_{G}(4D)\cdot (b_{8D})^{k}$ by the Bott isomorphism. $\text{Im}(i_{0}^{*})=(A+2D+6-2H)^{l}\cdot KO_{G}(4D)\cdot (b_{8D})^{k}$ for some $l\in \mathds{Z}_{\geq 0}$. A simple calculation shows that $\kappa o(X_{0})=2l$. Suppose $\kappa o(X_{1})=r$. Then we can find an element $z\in \widetilde{KO}_{G}(X_{1})$ such that $\varphi_{k_{1}}i_{1}^{*}(z)=2^{r}$. Therefore, $\varphi_{k_{0}}(\omega)=2^{r+\beta_{k_{1}}^{k_{1}-k_{0}}}$ where $\omega=(f^{S^{1}})^{*}(i_{1}^{*}(z))$. Since $k_{1}>k_{0}$, the map $(f^{S^{1}})^{*}$ factors through $KO_{G}((k_{0}+1)D)\rightarrow KO_{G}(k_{0}D)$. Therefore, we see that $\omega=\gamma(D)\cdot (a\eta(D)\lambda(D))\cdot (b_{8D})^{k}=a(1-D)\lambda(D)\cdot (b_{8D})^{k}$ for some $a\in \mathds{Z}$. Because of the commutative diagram, we have $\omega\in \text{Im}(i_{0}^{*})$. By Lemma \ref{lowest degree, odd}, we get $2^{2l+1}|\varphi_{k_{0}}(\omega)$. This implies $2l+1\leq r+\beta_{k_{1}}^{k_{1}-k_{0}}$.

(2) Suppose $X_{0}$ is even $KO_{G}$-split with $k_{0}=8k$. Notice that $\kappa o(X)=2l+2$ if $\mathcal{J}(X)$ is the submodule generated by $(2+A-2D-2B)^{l}(2-2D-B)(b_{8D})^{k}$ and $(A-2B)^{l}(B-A)(b_{8D})^{k}$. Using Lemma \ref{lowest degree, even}, the proof is almost the same with the previous case.
\end{proof}

   By Proposition \ref{periodicity of kappao}, we see that $\Sigma^{2H}X$ and $\Sigma^{8D}X$ are even (odd) $KO_{G}$-split if $X$ is even (odd) $KO_{G}$-split. Therefore, Proposition \ref{suspension by trivial} justifies the following definition:
\begin{defi}
A spectrum class $S=[(X,a,b+r)]$ with $a,b\in \mathds{Z},r\in[0,1)$ is called even (odd) $KO_{G}$-split if for integers $M,N\gg 0$, $\Sigma^{(8M-a)D}\Sigma^{(2N-b)H}X$ is even (odd) $KO_{G}$-split.
\end{defi}
\begin{ex}\label{sphere is split}
For any $a,b\in \mathds{Z}$ and $r\in[0,1)$, $[(S^{0},8a,2b+1+r)]$ is even $KO_{G}$-split and $[(S^{0},8a+4,2b+r)]$ is odd $KO_{G}$-split.
\end{ex}

The following proposition is easy to prove using Proposition \ref{split inequality for space}
\begin{pro}\label{split inequality for spectrum}
Let $S_{0},S_{1}\in \mathfrak{C}$ be two spectrum classes at level $k_{0},k_{1}$ respectively, with $k_{0}<k_{1}$. Suppose $S_{0}$ is even or odd $KO_{G}$-split and $S_{0}$ dominates $S_{1}$, then we have:
\begin{equation}
\kappa o(S_{0})<\kappa o(S_{1})+\beta_{k_{1}}^{k_{1}-k_{0}}.
\end{equation}
\end{pro}
Now let $Y$ be a homology sphere. Recall that we have a spectrum class invariant $S(Y)$ at level 0.
\begin{defi}\label{KO-split}
$Y$ is called Floer $KO_{G}$-split if $\Sigma^{H}S(Y)$ is even $KO_{G}$-split and $\Sigma^{4D}S(Y)$ is odd $KO_{G}$-split.
\end{defi}
\begin{rmk}
For simple examples like $Y=\pm\Sigma(2,3,12n+1) \text{ or }\pm\Sigma(2,3,12n+5)$, the two conditions in the above definition are either both true or both false. We expect that this fails in more complicated examples. If we only assume one of these two conditions, only half of the cases in Theorem \ref{split inequality for mfd} are still true.
\end{rmk}
\begin{rmk}\label{example of split mfd}
We will see in Section $8$ that $S^{3},\pm\Sigma(2,3,12n+1)$ and $-\Sigma(2,3,12n+5)$ are Floer $KO_{G}$-split, while $+\Sigma(2,3,12n+5)$ is not Floer $KO_{G}$-split.
\end{rmk}

\begin{proof}[Proof of Theorem \ref{split inequality for mfd}:]
(1) When $\mu(Y_{0})=0$, $S(Y_{0})=[(X,a,b)]$ for some space $X$ and some integers $a,b$. For large integers $M,N$, we have the following:

(i) The space $\Sigma^{(8M-a)D}\Sigma^{(2N-b+1)H}X$ is even $KO_{G}$-split.

(ii) The space $\Sigma^{(8M-a+4)D}\Sigma^{(2N-b)H}X$ is odd $KO_{G}$-split.

Now consider $p=4l+m$ for $m=0,1,2,3$:
\begin{itemize}
\item For $p=4l$, $\Sigma^{\frac{p}{2}H}\Sigma^{4D}S(Y_{0})=[(\Sigma^{4D}X,a,b-2l)]$ is odd $KO_{G}$-split by (ii).
\item For $p=4l+1$, $\Sigma^{\frac{p}{2}H}S(Y_{0})=[(\Sigma^{H}X,a,b-2l+\frac{1}{2})]$ is even $KO_{G}$-split by (i).
\item For $p=4l+2$,  $\Sigma^{\frac{p}{2}H}S(Y_{0})=[(\Sigma^{H}X,a,b-2l)]$ is even $KO_{G}$-split by (i).
\item For $p=4l+3$,  $\Sigma^{\frac{p}{2}H}\Sigma^{4D}S(Y_{0})=[(\Sigma^{4D}X,a,b-2l-2+\frac{1}{2})]$ is odd $KO_{G}$-split by (ii).

\end{itemize}
Similarly, we can prove that when $\mu(Y_{0})=1$, $\Sigma^{\frac{p}{2}H}S(Y_{0})$ is even $KO_{G}$-split for $p=4l+2$ and $4l+3$ while $\Sigma^{\frac{p}{2}H}\Sigma^{4D}S(Y_{0})$ is odd $KO_{G}$-split for $p=4l$ and $4l+1$.

Now repeat the proof of Theorem \ref{inequality} for $k=0$ or $4$, using Proposition \ref{split inequality for spectrum} instead of Proposition \ref{inequality for spectrum}. Notice that the two sides of the same inequalities are either both integers or both half-integers. The inequalities are proved.
\end{proof}

\section{Examples and Explicit bounds}

In this section, we will prove Theorem \ref{kappao for Brieskorn} about the values of $\kappa o_{i}(S^{3})$ and $\kappa o_{i}(\pm\Sigma(2,3,r))$ with $\text{gcd}(r,6)=1$. We will also use Corollary \ref{single boundary} to give some new bounds about the intersection forms of spin four manifolds with given boundaries.
\subsection{Basic Examples}
If $Y$ is a rational homology sphere admitting metric $g$ with a positive scaler curvature, then by the arguments in \cite{Manolescu1}, we obtain:
$$
S(Y,\mathfrak{s})=[(S^{0},0,n(Y,\mathfrak{s})/2)].
$$
In particular, $S^{3}$ is Floer $KO_{G}$-split and $\kappa o_{i}(S^{3})=0$ for any $i\in \mathds{Z}/8$.

In \cite{Manolescu2}, Manolescu gave two examples of spaces of type SWF that are related to the spectrum class invariants of the Brieskorn spheres $\pm\Sigma(2,3,r)$. We recall the construction here.

Suppose that $G$ acts freely on a finite $G$-CW complex $Z$, with the quotient space $Q=Z/G$. Let
$$\widetilde{Z}=([0,1]\times Z)/(0,z)\sim (0,z') \text{ and } (1,z)\sim (1,z') \text{ for all } z,z'\in Z$$
denote the unreduced suspension of $Z$, where $G$ acts trivially on the $[0,1]$ factor. We can take one of the two cone points (say $(0,z)\in \widetilde{Z}$) as the base point and view $\widetilde{Z}$ as a pointed $G$-space. It's easy to see that $\widetilde{Z}$ is of type SWF at level $0$.

We want to compute $\kappa o(\Sigma^{kD}\tilde{Z})$ for $k=0,1,...,7$. It turns out that the method in \cite{Manolescu2} also works here. Namely, the inclusion $i:(\Sigma^{kD}\tilde{Z})^{S^{1}}=\Sigma^{kD}S^{0} \rightarrow \Sigma^{kD}\tilde{Z}$ gives the long exact sequence:
\begin{equation}
...\rightarrow \widetilde{KO}_{G}(\Sigma^{kD}\tilde{Z})\stackrel{i^{*}}{\longrightarrow}  KO_{G}(kD)\stackrel{p^{*}}{\longrightarrow} KO_{G}^{1}(\Sigma^{kD}\tilde{Z},(kD)^{+})\rightarrow....
\end{equation}

By exactness of the sequence, we have $\text{Im}(i^{*})=\text{ker}(p^{*})$. By definition, we have: $$KO_{G}^{1}(\Sigma^{kD}\tilde{Z},(kD)^{+})\cong \widetilde{KO}_{G}^{1}(\Sigma^{kD}\Sigma Z_{+})\cong \widetilde{KO}_{G}(\Sigma^{kD}Z_{+}).$$ By abuse of notation, we still use $p^{*}$ to represent the map between $KO_{G}(kD)$ and $\widetilde{KO}_{G}(\Sigma^{kD}Z_{+})$. Checking the maps in the exact sequence, one can see that the $p^{*}$ is induced by the natural projection $p:\Sigma^{kD}Z_{+}\rightarrow (kD)^{+}$. Since $G$ acts freely on $\Sigma^{kD}Z_{+}$ away from the base point, we see that $\widetilde{KO}_{G}(\Sigma^{kD}Z_{+})\cong \widetilde{KO}((\Sigma^{kD}Z_{+})/G)$. Notice that $(Z\times kD)/G$ is a vector bundle over $Q$ and $(\Sigma^{kD}Z_{+})/G$ is the Thom space of this bundle. We are interested in two cases:
\begin{itemize}
\item $Z\cong G$, acting on itself via left multiplication.
\item $Z\cong T\cong S^{1}\times jS^{1}\subset \mathds{C}\times j\mathds{C}\subset \mathds{H}$ and $G$ acts on $T$ by left multiplication in $\mathds{H}$.
\end{itemize}

   The first case is easy since the isomorphism $\widetilde {KO}_{G}(\Sigma^{kD}Z_{+})\cong \widetilde {KO}(S^{k})$ is given by $i_{1}^{*}\circ r_{0}$, where $i_{1}: S^{k}\rightarrow \Sigma^{k\mathds{R}}Z_{+}$ is the standard inclusion and  $r_{0}:\widetilde {KO}_{G}(\Sigma^{kD}Z_{+})\rightarrow \widetilde {KO}(\Sigma^{k\mathds{R}}Z_{+})$ is the restriction map (See Fact \ref{restriction} in Section 2). It follows that $\text{Im}(i^{*})=\text{ker}(p^{*})=\text{ker}(i_{1}^{*}\circ r_{0}\circ p^{*})=\text{ker}(r)$, where $r:KO_{G}(kD)\rightarrow \widetilde{KO}(S^{k})$ is the restriction map.

We know the structure of $\widetilde{KO}(S^{k})$:
\begin{itemize}
\item $\widetilde{KO}(S^{0})\cong KO(pt)\cong \mathds{Z}$.
\item $\widetilde{KO}(S^{1})\cong \mathds{Z}_{2}$, generated by the Hurewicz image of the Hopf map in $\pi_{3}(S^{2})$.
\item $\widetilde{KO}(S^{2})\cong \mathds{Z}_{2}$, generated by the Hurewicz image of the square of the Hopf map.
\item $\widetilde{KO}(S^{4})\cong \mathds{Z}$, generated by $V_{\mathds{H}}-4$, where $V_{\mathds{H}}$ is the quaternion Hopf bundle.
\item $\widetilde{KO}(S^{k})\cong 0$ for $k=3,5,6,7$.
\end{itemize}

Therefore, by the explicit description of $\eta(D),\lambda(D),c$ after Theorem \ref{Pin(2)-gourp}. We get the following results about the kernel of $r:KO_{G}(kD)\rightarrow \widetilde{KO}(S^{k})$.
\begin{itemize}
\item For $k=0$, $\text{ker}(r)$ is the submodule generated by $1-D,A,B$.
\item For $k=1$, $\text{ker}(r)$ is generated by $2\eta(D)$.
\item For $k=2$, $\text{ker}(r)$ is generated by $2\eta(D)^{2}$ and $\gamma(D)^{2}c$.
\item For $k=4$, $\text{ker}(r)$ is  generated by $\lambda(D)-c, (1-D)\lambda(D), A\lambda(D)$ and $Ac$.
\item For $k=3,5,6,7$, $\text{ker}(r)\cong KO_{G}(kD)$.
\end{itemize}

From this, we get:
\begin{pro}\label{kappao for G}
$\kappa o(\Sigma^{kD}\tilde{G})=0$ for $k=3,4,5,6,7$ and $\kappa o(\Sigma^{kD}\tilde{G})=1$ for $k=0,1,2$.
\end{pro}

Now let's consider the case $Z\cong T$. We want to find $\text{ker}(p^{*})$ for $p^{*}:KO_{G}(kD)\rightarrow KO_{G}(\Sigma^{kD}T_{+})$. Notice that $S^{1}\subset G$ acts trivially on $(kD)^{+}$ and freely on $T$ with  $T/S^{1}=S^{1}$. We have $\widetilde{KO}_{G}(\Sigma^{kD}T_{+})=\widetilde{KO}((\Sigma^{kD}S^{1}_{+})/\mathds{Z}_{2})$. The space $(\Sigma^{kD}S^{1}_{+})/\mathds{Z}_{2}$ can be identified with:
$$
[0,1]\times(kD)^{+}/(0,x)\sim(1,-x) \text{ and }(t_{1},\infty)\sim(t_{2},\infty) \text{ for any } x\in (kD)^{+} \text{ and } t_{1},t_{2}\in[0,1].
$$

Consider the inclusion $i_{2}:\{0\}\times (kD)^{+}\rightarrow (\Sigma^{kD}S^{1}_{+})/\mathds{Z}_{2}$. Notice that $((\Sigma^{kD}S^{1}_{+})/\mathds{Z}_{2})/(kD)^{+}$ $\cong S^{k+1}$. We get the long exact sequence:
\begin{equation}\label{exact sequence}...\rightarrow \widetilde{KO}(S^{k+1}) \stackrel{\delta}{\longrightarrow}\widetilde{KO}(S^{k+1})\rightarrow \widetilde{KO}((\Sigma^{kD}S^{1}_{+})/\mathds{Z}_{2})\stackrel{i_{2}^{*}}{\longrightarrow}\widetilde{KO}(S^{k})\rightarrow...\end{equation}

By checking the iterated mapping cone construction, which gives us this long exact sequence, it is not hard to prove that $\delta$ is induced by the map $f:S^{k+1}\rightarrow S^{k+1}$ with $\text{deg}(f)=0$ for even $k$ and $\text{deg}(f)=2$ for odd $k$.

When $k=2,4,5,6$, we have $\widetilde{KO}(S^{k+1})=0$. Therefore, $i_{2}^{*}$ is injective, which implies $i_{1}^{*}\circ r_{0}:\widetilde{KO}_{G}(\Sigma^{kD}T_{+})\rightarrow \widetilde{KO}((kD)^{+})$ is injective ($i_{1}^{*}$ and $r_{0}$ are defined as in the case $Z\cong G$). We see that when $k=2,4,5,6$, just like the case $Z\cong G$, the kernel of $p^{*}$ is the kernel of the restriction map $r:KO_{G}(kD)\rightarrow \widetilde{KO}(S^{k})$. Thus, we get $\kappa o(\Sigma^{kD}\tilde{T})=\kappa o(\Sigma^{kD}\tilde{G})$ for $k=2,4,5,6$.

For $k=0$, consider $[0,1]$ as the subset $\{1+je^{i\theta}|\theta\in [0,\pi]\}\subset T$. The left endpoint is mapped to the right endpoint under the action of $-j\in G$. This embedding of $[0,1]$ gives us the following explicit description of the map $p^{*}:RO(G)\cong \widetilde{KO}_{G}(S^{0})\rightarrow \widetilde{KO}_{G}(T_{+})\cong KO_{G}(T)\cong KO(T/G)=KO(S^{1})$.

Starting from a representation space $V$ of $G$, we get an trivial bundle $V\times [0,1]$ over $[0,1]$. Identifying $(x,0)$ with $((-j)\circ x,1)$ for any $x\in V$, we get a bundle $E$ over $S^{1}$. $[E]\in KO(S^{1})$ is the image of $[V]\in RO(G)$ under $p^{*}$.

We know that $KO(S^{1})$ is generated by  the one dimensional trivial bundle [1] and the one dimensional nontrivial bundle $[m]$, subject to the relation 2([1]-[m])=0. Using the explicit description of $p^{*}$, we see that $p^{*}(1)=[1]$, $p^{*}(D)=[m]$ and $p^{*}(A)=p^{*}(B)=0$. Therefore, we get $\kappa o(\tilde{T})=2$.

Applying Corollary \ref{relation} for $S=\Sigma^{2D}\tilde {T}$, we get $\kappa o(\Sigma^{3D}\tilde{T})+1\geq \kappa o(\Sigma^{2D}\tilde{T})=1$. Applying Corollary \ref{relation} for $S=\Sigma^{3D}\tilde {T}$, we get $0=\kappa o(\Sigma^{4D}\tilde{T})+0\geq \kappa o(\Sigma^{3D}\tilde{T})$. Therefore, we see that $\kappa o(\Sigma^{3D}\tilde{T})=0$.

Applying Corollary \ref{relation} for $S=\Sigma^{2D}\tilde{T}$ and $S=\Sigma^{D}\tilde{T}$, we get $\kappa o(\Sigma^{D}\tilde{T})=1$ or $2$.

 For $k=7$, the map $\delta:\widetilde{KO}(S^{8})\rightarrow \widetilde{KO}(S^{8})$ is multiplication by $2$. Since $\widetilde{KO}(S^{7})=0$, we get $\widetilde{KO}((\Sigma^{kD}S^{1}_{+})/\mathds{Z}_{2})=\mathds{Z}_{2}$. This implies $p^{*}(2b_{8D}\cdot\gamma(D))=2p^{*}(b_{8D}\cdot\gamma(D))=0$. Therefore, $2b_{8D}\cdot\gamma(D)\in \text{ker}(p^{*})$ and $\kappa o(\Sigma^{7D}\tilde{T})=0$ or $1$.

\begin{lem}
$\kappa o(\Sigma^{D}\tilde{T})=2$ and $\kappa o(\Sigma^{7D}\tilde{T})=1$.
\end{lem}

\begin{proof}
This can be proved directly using Gysin sequence. But here we use a different approach. In \cite{Manolescu2} and \cite{Manolescu4}, Manolescu proved that $S(-\Sigma(2,3,11))=[(\tilde{T},0,1)]$, where $-\Sigma(2,3,11)$ is a negative oriented Brieskorn sphere. Therefore, by Definition \ref{kappao for spectrum} and Proposition \ref{periodicity of kappao for spectrum}, we get: $$\kappa o_{i}(-\Sigma(2,3,11))=\kappa o(\Sigma^{(i+4)D}\tilde{T})+1-\beta_{i}^{4}. $$ In particular, $\kappa o_{3}(-\Sigma(2,3,11))=\kappa o(\Sigma^{7D}\tilde{T})-2$ and $\kappa o_{5}(-\Sigma(2,3,11))=\kappa o(\Sigma^{D}\tilde{T})-2$. Since $-\Sigma(2,3,11)$ bounds a smooth spin four manifold with intersection from $\left(\begin{smallmatrix}
 0 & 1 \\
 1& 0
\end{smallmatrix}\right)$ (see \cite{Manolescu2}). We can apply Corollary \ref{single boundary} for $p=0,q=1$ and get $\kappa o_{5}(-\Sigma(2,3,11))\geq 0$, which implies $\kappa o(\Sigma^{D}\tilde{T})\geq 2$. We get $\kappa o(\Sigma^{D}\tilde{T})=2$ by our discussion before the lemma.

We can also apply Theorem \ref{inequality} for $Y_{0}=S^{3},Y_{1}=-\Sigma(2,3,11),p=0,q=1$ and $k=2$. We have $\kappa o_{3}(-\Sigma(2,3,11))\geq -1$ and $\kappa o(\Sigma^{7D}\tilde{T})\geq 1$. Therefore, $\kappa o(\Sigma^{7D}\tilde{T})=1$ by our discussions before.
\end{proof}
We summarise our results in the following proposition.
\begin{pro}\label{kappao for T}
$\kappa o(\Sigma^{kD}\tilde{T})=2$ for $k=0,1$; $\kappa o(\Sigma^{kD}\tilde{T})=1$ for $k=2,7$ and $\kappa o(\Sigma^{kD}\tilde{T})=0$ for $k=3,4,5,6$.
\end{pro}

Now we calculate $\kappa o_{i}(\pm \Sigma(2,3,r))$ with $\text{gcd}(6,r)=1$. Actually, the spectrum class invariants $S(\pm\Sigma(2,3,r))$ are given in \cite{Manolescu2}.

\begin{pro}[Manolescu \cite{Manolescu2}]We have the following results about $S(\pm \Sigma(2,3,r))$.
 \begin{itemize}
 \item$S(\Sigma(2,3,12n-1))=[(\tilde{G}\vee\underbrace{\Sigma G_{+}\vee...\vee\Sigma G_{+}}_{n-1},0,0)].$
 \item  $S(-\Sigma(2,3,12n-1))=[(\tilde{T}\vee\underbrace{\Sigma^{2} G_{+}\vee...\vee\Sigma^{2} G_{+}}_{n-1},0,1)].$
\item $S(\Sigma(2,3,12n-5))=[(\tilde{G}\vee\underbrace{\Sigma G_{+}\vee...\vee\Sigma G_{+}}_{n-1},0,1/2)].$
\item $S(-\Sigma(2,3,12n-5))=[(\tilde{T}\vee\underbrace{\Sigma^{2} G_{+}\vee...\vee\Sigma^{2} G_{+}}_{n-1},0,1/2)].$
\item$S(\Sigma(2,3,12n+1))=[(S^{0}\vee\underbrace{\Sigma^{-1} G_{+}\vee...\vee\Sigma^{-1} G_{+}}_{n},0,0)].$\footnote{Strictly speaking, by this we mean the spectrum class of $(\mathds{H}^{+}\vee\underbrace{\Sigma^{3} G_{+}\vee...\vee\Sigma^{3} G_{+}}_{n},0,1)$.}
\item$S(-\Sigma(2,3,12n+1))=[(S^{0}\vee\underbrace{G_{+}\vee...\vee G_{+}}_{n},0,0)].$
\item$S(\Sigma(2,3,12n+5))=[(S^{0}\vee\underbrace{\Sigma^{-1}G_{+}\vee...\vee\Sigma^{-1}G_{+}}_{n},0,-1/2)].$
\item$S(-\Sigma(2,3,12n+5))=[(S^{0}\vee\underbrace{G_{+}\vee...\vee G_{+}}_{n},0,1/2)].$
\end{itemize}
\end{pro}

As we mentioned in Remark \ref{example of split mfd}, $\pm\Sigma(2,3,12n+1)$ and $-\Sigma(2,3,12n+5)$ are $KO_{G}$-split because of Example \ref{sphere is split}. While using the relations in Theorem \ref{Pin(2)-gourp} and Theorem \ref{multiplicative structure}, it is not hard to prove that the space $(8MD\oplus(2N+2)H)^{+}$ is not even $KO_{G}$-split for integers $M,N\gg0$. This implies that $+\Sigma(2,3,12n+5)$ is not $KO_{G}$-split.

Since it's easy to see that wedging with a free $G$-space does not change the $\kappa o$ invariants, we don't need to consider those $\Sigma^{l} G_{+}$ factors. By Definition \ref{kappao for spectrum} and Proposition \ref{periodicity of kappao for spectrum}, we can use Proposition $\ref{kappao for G}$ and Proposition $\ref{kappao for T}$ to prove the results in Theorem \ref{kappao for Brieskorn} easily.
\subsection{Explicit Bounds}
Now we use Corollary \ref{single boundary} and Proposition \ref{stable map between spheres} to get explicit bounds on the intersection forms of spin $4$-manifolds with boundary $\pm\Sigma(2,3,r)$.
\begin{thm}\label{explicit bounds}
Let $W$ be an oriented, smooth spin $4$-manifold with $\partial W=\pm \Sigma(2,3,r)$. Assume that the intersection form of $W$ is $p(-E_{8})\oplus q\left(\begin{smallmatrix}
 0 & 1 \\
 1& 0
\end{smallmatrix}\right)$ for $p>1,q>0$.\footnote{It is easy to see that the conclusions are not true for $p=0,1$. For example, $\pm\Sigma(2,3,12n-1)$ bounds a spin manifold with intersection form $\left(\begin{smallmatrix}
 0 & 1 \\
 1& 0
\end{smallmatrix}\right)$.} If the mod $8$ reduction of $p$ is $m$, then we have $q-p\geq c_{m}$, where $c_{m}$ are constants listed below. (Recall that the mod $2$ reduction of $p$ is the Rohklin invariant of the boundary.)

\begin{center}
    \begin{tabular}{ |c|c|c |c |c|}
     \hline& $m=0$ & $m=2$ & $m=4$ & $m=6$ \\ \hline
    $\Sigma(2,3,12n-1)$ & $2$ & $0$ & $1$ & $2$ \\ \hline
    $-\Sigma(2,3,12n-1)$ & $3$ & $(2)$ & $(3)$ & $3$ \\ \hline
    $\Sigma(2,3,12n+1)$ & $(3)$ & $1$ & $(2)$ & $(3)$ \\ \hline
    $-\Sigma(2,3,12n+1)$ & $3$ & $1$ & $2$ & $3$ \\ \hline
    \end{tabular}
\end{center}
\begin{center}
    \begin{tabular}{ |c|c|c |c |c|}
     \hline& $m=1$ & $m=3$ & $m=5$ & $m=7$ \\ \hline
    $\Sigma(2,3,12n-5)$ & $1$ & $2$ & $3$ & $3$ \\ \hline
    $-\Sigma(2,3,12n-5)$ & $2$ & $(1)$ & $(2)$ & $2$ \\ \hline
    $\Sigma(2,3,12n+5)$ & $(2)$ & $0$ & $(1)$ & $(2)$ \\ \hline
    $-\Sigma(2,3,12n+5)$ & $2$ & $3$ & $4$ & $4$ \\ \hline
    \end{tabular}
\end{center}
\end{thm}
\begin{rmk}
Some of the bounds in Theorem \ref{explicit bounds} can also be obtained by other methods. For example, the case $m=2$ for $\Sigma(2,3,12n+1)$ can be obtained using $\kappa$-invariant (see \cite{Manolescu2}). Also, some bounds can be obtained by filling method for small $n$. For example, the case $m=2,4$ for $-\Sigma(2,3,11)$ can be deduced from Theorem \ref{inequality for closed mfd}, using the fact that $\Sigma(2,3,11)$ bounds a spin $4$-manifold with intersection form $2(-E_{8})\oplus 2\left(\begin{smallmatrix}
 0 & 1 \\
 1& 0
\end{smallmatrix}\right)$. However, the bounds that we put in the brackets in Theorem \ref{explicit bounds} appear to be new for general $n$.
\end{rmk}

\begin{proof}
Since we can do surgeries on loops without changing intersection forms, we will always assume $b_{1}(W)=0$.

 (1) Suppose $\Sigma(2,3,12n+1)$ bounds a spin $4$-manifold with intersection form $8l(-E_{8})\oplus (8l+2)\left(\begin{smallmatrix}
 0 & 1 \\
 1& 0
\end{smallmatrix}\right)$ for $l>0$. Then we get a spin cobordism from $-\Sigma(2,3,12n+1)$ to $S^{3}$ with the same intersection form. By Theorem \ref{cobordism map}, $\Sigma^{4lH}S(-\Sigma(2,3,12n+1))$ dominates $\Sigma^{8l+2}S(S^{3})$.
Since $S(-\Sigma(2,3,12n+1))=[(S^{0}\vee G_{+}\vee...\vee G_{+},0,0)]$ and $S(S^{3})=[(S^{0},0,0)]$, by Definition \ref{stable equivalence}, we get a map: $$f:\Sigma^{r\mathds{R}+(4l+M)H+ND}(S^{0}\vee G_{+}\vee...\vee G_{+})\rightarrow \Sigma^{r\mathds{R}+MH+(8l+2+N)D}S^{0}$$ for some $M,N\in \mathds{Z}$. 
Restricting to the first factor of $S^{0}\vee G_{+}\vee...\vee G_{+}$, we obtain:
$$g:\Sigma^{r\mathds{R}+(4l+M)H+ND}S^{0}\rightarrow \Sigma^{r\mathds{R}+MH+(8l+2+N)D}S^{0},$$
which induces homotopy equivalence between the $G$-fixed point sets. This a contradiction with Proposition \ref{stable map between spheres}. The case $m=0$ for $\Sigma(2,3,12n+1)$ is proved.

(2) Suppose $\Sigma(2,3,12n+5)$ bounds a smooth spin manifold with intersection form $(8l+1)(-E_{8})\oplus (8l+2)\left(\begin{smallmatrix}
 0 & 1 \\
 1& 0
\end{smallmatrix}\right)$ for $l>0$. Then we get a spin cobordism from $-\Sigma(2,3,12n+5)$ to $S^{3}$. As the previous case, this implies $\Sigma^{(4l+1/2)H}S(-\Sigma(2,3,12n+5))$ dominates $\Sigma^{(8l+2)D}S(S^{3})$. Since $\Sigma^{(4l+1/2)H}S(-\Sigma(2,3,12n+5))=[(\Sigma^{4lH}S^{0},0,0)]$, we get the contradiction as before. This proves the case $m=1$ for $\Sigma(2,3,12n+5)$.

 (3) Suppose $-\Sigma(2,3,12n-1)$ bounds a spin $4$-manifold with intersection form $(8l+2)(-E_{8})\oplus(8l+3)\left(\begin{smallmatrix}
 0 & 1 \\
 1& 0
\end{smallmatrix}\right)$ for $l\geq0$. By Corollary \ref{single boundary}, we get $4l+3<\kappa o_{3+8l}(-\Sigma(2,3,12n-1))+\beta_{8l+3}^{8l+7}=-1+4+4l$, which is a contradiction. This proves the case $m=2$ for $-\Sigma(2,3,12n-1)$.

Using similar method as (3), we can prove all the other cases except:
\begin{itemize}
\item $m=0$ for $\pm \Sigma(2,3,12n-1)$ and $-\Sigma(2,3,12n+1)$,
 \item $m=7$ for $\Sigma(2,3,12n-5)$ and $-\Sigma(2,3,12n+5)$,
 \item  $m=1$ for $-\Sigma(2,3,12n-5)$.
\end{itemize}

(4) We need to introduce another approach in order to prove the rest of the cases. Consider the orbifold $D^{2}$-bundle over $S^{2}(2,3,r)$. This gives us an orbifold $X'$ with boundary $+\Sigma(2,3,r)$. We have $b_{2}^{+}(X')=0,b_{2}^{-}(X)=1$ and $X'$ has a unique spin structure $\mathfrak{t}$. Now suppose $-\Sigma(2,3,r)$ bounds a spin manifold $X$ with intersection form $p(-E_{8})\oplus q\left(\begin{smallmatrix}
 0 & 1 \\
 1& 0
\end{smallmatrix}\right)$. Then we can glue $X$ and $X'$ together to get an oriented closed spin $4$-orbifold. We have:
$$ \text{ind}_{\mathds{C}}\slashed{D}(X\cup X')=p+\omega(\Sigma(2,3,r),X',\mathfrak{t}).$$
Here $\omega(\Sigma(2,3,r),X',\mathfrak{t})$ is the Fukumoto-Furuta invariant defined in \cite{Fukumoto-Furuta}. Saveliev \cite{Saveliev} proved that $\omega(\Sigma(2,3,r),X',\mathfrak{t})=-\overline{\mu}(\Sigma(2,3,r))=\overline{\mu}(-\Sigma(2,3,r))$, where $\overline{\mu}$ is the Neumann-Siebenmann invariant \cite{Neumann-Siebenmann1,Neumann-Siebenmann3}. In \cite{Fukumoto-Furuta}, Fukumoto and Furuta considered the finite dimensional approximation of the Seiberg-Witten equations on the orbifold $X\cup X'$ and constructed a stable Pin($2$)-equivariant map: $(\frac{\text{ind}_{\mathds{C}}\slashed{D}(X\cup X')}{2} H)^{+}\rightarrow (b_{2}^{+}(X\cup X') D)^{+}$ which induces homotopy equivalence on the Pin($2$)-fixed point set. (Recall that $H$ and $D$ are Pin($2$)-representations defined in Section 2). Since $b_{2}^{+}(X\cup X')=q$ and $\text{ind}_{\mathds{C}}\slashed{D}(X\cup X')=p+\overline{\mu}(-\Sigma(2,3,r))$, we can apply Proposition \ref{stable map between spheres} to get:
 $$q-p\geq 3+\overline{\mu}(-\Sigma(2,3,r)) \text{ if } 0<p+\overline{\mu}(-\Sigma(2,3,r)) \text{ can be divided by }8.$$

Similarly, suppose $\Sigma(2,3,r)$ bounds a spin $4$-manifold $X'$ with intersection form $p(-E_{8})\oplus q\left(\begin{smallmatrix}
 0 & 1 \\
 1& 0
\end{smallmatrix}\right)$. We can consider $X'\cup (-X)$ and repeat the argument above. We get:
$$q-p\geq 2+\overline{\mu}(\Sigma(2,3,r)) \text{ if } 0<p+\overline{\mu}(\Sigma(2,3,r)) \text{ can be divided by }8.$$

The invariants $\overline{\mu}(\pm\Sigma(2,3,r))$ were computed in \cite{Neumann-Siebenmann1,Neumann-Siebenmann3}: $$\overline{\mu}(\pm\Sigma(2,3,12n-1))=\overline{\mu}(\pm\Sigma(2,3,12n+1))=0,$$  $$\overline{\mu}(\Sigma(2,3,12n-5))=\overline{\mu}(-\Sigma(2,3,12n+5)=1,$$
$$\overline{\mu}(-\Sigma(2,3,12n-5))=\overline{\mu}(\Sigma(2,3,12n+5)=-1.$$
Therefore, simple calculations prove the rest of the cases.
\end{proof}






\vskip 0.5 truecm

\end{document}